\documentclass[12pt,leqno]{amsart}
\usepackage{amsmath}
\usepackage{amssymb}
\usepackage{amsthm} 
\usepackage{epsf}
\usepackage{graphicx}
\usepackage{esint} 
\usepackage{dsfont}
\usepackage[pagebackref]{hyperref} 
\hypersetup{backref,  pdfpagemode=FullScreen,  colorlinks=true} 
\usepackage{color}
\usepackage{enumerate}

\numberwithin{equation}{section}

\newcommand{\nn}{\nonumber}
\newcommand{\ms}{\medskip}

\newcommand{\R}{\mathbb{R}}
\renewcommand{\H}{\mathcal H}

\renewcommand{\d}{\partial}
\newcommand{\dist}{\,\mathrm{dist}}

\newcommand{\sm}{\setminus}
\DeclareMathOperator{\supp}{supp}
\DeclareMathOperator{\diam}{diam}

\newcommand{\ol}{\overline}
\newcommand{\ub}{\underbar}

\newcommand{\cF}{{\mathcal  F}}

\newcommand{\cG}{{\mathcal  G}}

\newcommand{\cB}{{\mathcal  B}}

\newcommand{\1}{{\mathds 1}}

\theoremstyle{plain}
\newtheorem{theorem}[equation]{Theorem}
\newtheorem{lemma}[equation]{Lemma}

\newtheorem{proposition}[equation]{Proposition}

\theoremstyle{definition}
\newtheorem{definition}[equation]{Definition}

\theoremstyle{remark}
\newtheorem{remark}[equation]{Remark}

\newcommand{\dv}{\operatorname{div}}

\newcommand{\RR}{{\mathbb{R}}}

\newcommand{\vp}{\varphi}

\newcommand{\po}{{\partial\Omega}}

\newcommand{\bp}{\noindent {\it Proof}.\,\,}
\newcommand{\ep}{\hfill$\Box$ \vskip 0.08in}

\def\div{\mathop{\operatorname{div}}}

\definecolor{greenish}{cmyk}{0.91,0,0.88,0.12}

\newcounter{hcountno}

\begin{document}
~
\title[UR and approximation of Green functions by distance functions]{Approximation of Green functions and domains with uniformly rectifiable boundaries
of all dimensions}

\author[David]{Guy David}
\address{Guy David. Universit\'e Paris-Saclay, CNRS, Laboratoire de math\'ematiques d'Orsay, 91405 Orsay, France}
\email{guy.david@universite-paris-saclay.fr}

\author[Mayboroda]{Svitlana Mayboroda}
\address{Svitlana Mayboroda. School of Mathematics, University of Minnesota, Minneapolis, MN 55455, USA}
\email{svitlana@math.umn.edu}

\thanks{David was partially supported by the European Community H2020 grant GHAIA 777822, and the Simons Foundation grant 601941, GD. Mayboroda was supported in part 
by the NSF grant DMS 1839077 and the Simons foundation grant 563916, SM}
\date{} 

\maketitle

\begin{abstract}
The present paper establishes equivalence between uniform rectifiability of 
the boundary of a   domain and the property that the Green function for elliptic operators is well approximated by affine functions (distance to the hyperplanes). 
The results are novel in a variety of ways, in particular
\begin{enumerate}
\item this is the first time the underlying property of the control of the Green function by affine functions, or by the distance to the boundary, in the sense of  the Carleson prevalent sets, appears in the literature; the ``direct" result  established here is new even in the half space;
\item the results are optimal, providing a full characterization of uniform rectifiability  under the (standard) mild topological assumptions;
\item to the best of the authors' knowledge, this is the first free boundary result applying to all elliptic operators, without any restriction on the coefficients 
(the direct one assumes the standard, and necessary, Carleson measure condition); 
\item our theorems apply to all domains, with possibly lower dimensional boundaries: this is the first free boundary result in higher co-dimensional setting and as such, the first PDE characterization of uniform rectifiability for a set of dimension $d$, $d<n-1$, in $\RR^n$.
\end{enumerate}
The paper offers a general way to deal with related issues
considerably beyond the scope of the aforementioned theorem, including the 
question  of approximability of the {\it gradient} of the Green function, and the comparison of the Green function to a certain version of the distance to the original set rather than distance to the hyperplanes.  
\end{abstract}

\ms\noindent{\bf Key words:}  
Green function, Elliptic operators, Uniform rectifiability
\ms\noindent

\tableofcontents

\section{Introduction}
We consider elliptic operators, or their generalizations, on a domain $\Omega \subset \R^n$.
The main goal of the present paper is to show that $\Omega$ 
is reasonably regular (uniformly rectifiable) if and only if the Green function can be well approximated by distances to planes, or by a certain distance to 
the boundary $\d\Omega$.
Formally speaking, the result is one of the end-points of the big quest of establishing sharp, 
optimal connections between the geometric and PDE properties of sets -- see the work by J.~Garnett, S.~Hofmann, J.-M.~Martell, X.~Tolsa, T.~Toro, and their collaborators 
\cite{HMU, HMMTZ, AHMMT}. However, it stands apart both philosophically and in terms of the involved techniques, and it makes 
a significant 
breakthrough in several directions. First of all, perhaps surprisingly, the underlying property of 
good approximation of the Green function by 
distance functions has never appeared in the literature{\footnote{In fact, the idea that the Green function is almost affine whenever the set is reasonably flat permeates virtually all  approaches to free boundary problems, including the work of Alt and Caffarelli, Kenig and Toro, Lewis and Vogel, Hofmann, Martell, Tolsa, Garnett, Nystrom, and others. Yet, the particular condition using Carleson prevalent sets which gives equivalence with uniform rectifiability is new as far as the authors know.}}. The ``direct" result states that for a certain class 
of elliptic operators, the Green function is morally affine in the sense of Carleson prevalent sets. This  is new even in a half space, not to mention in a domain with uniformly rectifiable boundary. Secondly, the ``free boundary" result, stating that the desired approximability of the Green function  by distances to planes implies uniform rectifiability  is, to the best  of our knowledge, the first free boundary result in this context that pertains to all elliptic operators, without any additional restrictions on the coefficients of the equation. Finally, and secretly this was our main motivation, this paper is the first PDE characterization of the lower-dimensional uniformly rectifiable sets. 

The efforts  in this direction have spanned now about a decade and could roughly be split  
into two principal directions: the work on reflectionless measures by F. Nazarov, B. Jaye, 
and their collaborators, which unfortunately still hinges on the problem of proving the reflectionless property for the key operators at hand, 
and the work by the authors of the present paper and J. Feneuil which identifies 
an appropriate PDE context and proves a number of ``direct" results,  but also faces a mysterious problem in the ``free boundary" direction. Indeed, we discovered that in  domains with lower dimensional boundaries there exist elliptic operators for which the elliptic measure is proportional to the Hausdorff measure on any AR set, independently of any sort of regularity  or uniform  rectifiability, which  shatters the natural conjectures inspired by the co-dimension one case. This is what forced us to turn to the Green function and its proximity to the affine functions, or to the distance: such a condition seems to control more accurately the torsion of the solution around a low-dimensional boundary than the more traditional estimates on the harmonic measure. It appears, for all the reasons outlined above, that in a certain sense it is stronger and more appropriate in the classical domains with $n-1$ dimensional boundary as well. Let us now turn to the details. 

In all this paper, $\Omega$ will be a domain in $\R^n$ whose boundary 
$E = \d \Omega$ is Ahlfors regular of some dimension $d < n$. 
This means that there is a measure $\mu$, supported
on $E$, and a constant $C \geq 1$ such that 
\begin{equation} \label{1a1}
C^{-1} r^d \leq \mu(B(x,r)) \leq C r^d  
\ \text{ for $x\in E$ and $0 < r < +\infty$.} 
\end{equation}
It is well known that in this case $\mu$ is equivalent to $\H^d_{\vert E}$, 
the restriction to $E$ of the Hausdorff measure of dimension $d$. 
In this paper we assume that $\Omega$ and $E$ are unbounded; the case of bounded 
domains would be similar but need a little more notation. 
Note that a priori we do not assume 
that $d$ is integer and allow, 
for instance, Ahlfors regular snowflakes and Cantor sets. Instead, we will prove that the desired approximation of the Green function (when meaningful) implies that the dimension is integer and the set is rectifiable. 

When $d \geq n-1$, we shall systematically add the assumption that $\Omega$ is a one-sided $NTA$ domain, which means that we assume some quantitative openness and connectedness in the form of  the existence of corkscrew points  and Harnack chains in $\Omega$; see the beginning of Section \ref{S2} for the definitions. 
This is a traditional topological background hypothesis in this context. Fortunately, 
we do not need to impose this explicitly when $d<n-1$, 
because then $\Omega$  satisfies
the corkscrew and Harnack chain conditions automatically. 

Turning to the operators, we will split the discussion into two cases. 
When $n-2<d<n$, we will concentrate  on the classical Laplacian $\Delta$, or more generally 
elliptic operators $L = - \dv A \nabla$, with bounded measurable, not necessarily symmetric, coefficients -- see \eqref{2b2}--\eqref{2b3} for the (usual) definition of ellipticity. 

When $d\leq n-2$, the classical elliptic operators are not appropriate. Their solutions do not ``see" the lower dimensional sets, and for instance, a bounded  harmonic function in $\RR^n\setminus \RR^d$ is indistinguishable from a harmonic 
function in $\RR^n$. Over the recent years, the authors of the present paper, 
together with J. Feneuil, M. Engelstein, and other collaborators, developed a rather complete 
elliptic theory on such domains which identifies a certain class of degenerate elliptic operators 
as a proper substitute to 
standard elliptic operators
in this setting \cite{DFM1, DFM2, DFM3, DFM5, DM1, Fen}. 
It was shown that the general results, such as the maximum principle, 
trace and extension theorems, existence of the harmonic measure and Green function, 
all hold for the operators 
$$L=-\dv A \dist (\cdot, E)^{d+1-n} \nabla, $$
where $A$ is the usual elliptic matrix as above and $\dist (\cdot, E)$ is the Euclidean distance to the boundary. Note that for $d=n-1$ these  
operators coincide with the classical elliptic ones. 

It turns out, however, perhaps surprisingly, that the analogue of the emblematic ``Laplacian" 
case is not simply $L=-\dv  \dist (\cdot, E)^{d+1-n} \nabla,$ unless $d=n-1$, 
in the sense that the Euclidean distance seems 
too rough to establish, for instance, the 
solvability of boundary problems on non-smooth domains, or the 
absolute continuity of harmonic measure with respect to the Lebesgue measure, and one has to be more creative than that. 

To this end, our favorite degenerate elliptic operator has become
\begin{equation} \label{1a2}
L = L_\alpha = - \dv D_\alpha^{d+1-n} \nabla,
\end{equation}
where $\alpha > 0$ is a parameter and $D_\alpha$ is a smoother distance function defined by
\begin{equation} \label{1a3}
D_\alpha(X) = D_{\alpha,\mu}(X)  = R_{\alpha,\mu}(X)^{-1/\alpha},
\end{equation}
with
\begin{equation} \label{1a4}
R_{\alpha,\mu}(X)  = \int_E |X-y|^{-d-\alpha} d\mu(y),
\end{equation}
and where $\mu$ is any Ahlfors regular (AR for short) measure on $E = \d \Omega$, i.e., 
any measure that satisfies \eqref{1a1}. 
This is the best substitute for the Laplacian that we could find; 
in particular, we have showed that for such an operator the elliptic measure 
is an $A^\infty$ weight for any domain with a uniformly rectifiable boundary \cite{DM1, Fen}.

As we have mentioned above, our free boundary results will hold in the general setting of elliptic operators as above. The ``direct" results are, by necessity, restricted to operators which are morally similar to the Laplacian or, in the case of lower dimensional boundaries, to \eqref{1a2}. In either setting, the resulting class of operators is of the nature of the best possible. In fact, it is slightly more general than anything previously considered in this context, but morally speaking we impose the ``usual" Carleson measure condition on the coefficients, whose failure is known to produce abundant counterexamples \cite{MM, P}, see also Remark \ref{r3b1}. To state this more precisely, we need some definitions.

We shall focus in this paper on so-called weak properties, whose model is always the same. 
We choose some desirable property (either geometric or 
related to the Green function), often stated in terms of some parameters like $\varepsilon > 0$, 
then consider the set of balls $B(x,r)$ for which it fails, and require that this 
bad set satisfies a Carleson packing condition (we will also say that the complement 
is a Carleson-prevalent set), as follows.

\begin{definition} \label{d1}
Let $\cB$ be a subset of $E \times (0,+\infty)$. We say that $\cB$ satisfies a 
\ub{Carleson} \ub{packing condition}
when there is a constant $C \geq 0$ such that for every $x \in E$ and $r > 0$,
\begin{equation} \label{1a6}
\int_{y \in E \cap B(x,r)}\int_{0<t<r} \1_{\cB}(y,t) \frac{d\mu(y) dt}{t} \leq C r^d.
\end{equation}
We say that $\cG \subset E \times (0,+\infty)$ is a \ub{Carleson-prevalent set} when 
$E \times (0,+\infty)\sm \cG$ satisfies a Carleson packing condition.
\end{definition} 

The term prevalent is new in this context (we think it will help to give a name to the good 
sets too), and has no intended relation with existing uses of prevalence in other domains 
of mathematics. Here we chose an Ahlfors regular measure $\mu$ on $E$, but if we 
chose a different one, for instance $\H^d_{\vert E}$, we obtained the same notion 
(but a different Carleson constant $C$).

\ms

Given an elliptic operator $L = - \dv A \nabla$ we say  that $L$ is sufficiently close 
locally to a constant coefficient elliptic operator, if the following weak Carleson measure 
condition holds. 
For each choice of constants $\tau > 0$ (small) and 
$K \geq 1$ (large), denote by $\cG_{cc}(\tau,K)$ the set of pairs
$(x,r)\in E \times (0,+\infty)$ such that there is a constant matrix $A_0 = A_0(x,r)$
such that 
\begin{equation} \label{3b1}
\int_{X \in W_K(x,r)} |A(X) - A_0| dX \leq \tau r^n,
\end{equation}
where $W_K(x,r)$ is a large Whitney region associated to $B(x,r)$, defined by
\begin{equation} \label{3b2}
W_K(x,r) = \big\{ X \in \Omega \cap B(x,Kr) \, ; \, \dist(X,E) \geq K^{-1} r \big\}.
\end{equation}
Our condition is that
\begin{equation} \label{3b3}
\text{for every choice of $\tau > 0$ and $K \geq 1$, $\cG_{cc}(\tau,K)$ is a 
Carleson prevalent set}
\end{equation}
(as in Definition \ref{d1}), or in other terms the corresponding bad set 
$\cB_{cc}(\tau,K) = E \times (0,+\infty) \sm \cG_{cc}(\tau,K)$
satisfies a Carleson packing condition.

Incidentally, we can assume that $A_0$ is elliptic with the same ellipticity constant  
as the $A(X)$, because if \eqref{3b1} holds for any $A_0$, it also holds (with a worse constant) 
with $A_0 = A(X_0)$, for some $X_0 \in W_K$ chosen by Chebyshev.

It would be easy to check that this condition is weaker than the standard conditions used 
in similar circumstances, such as the Dahlberg-Kenig-Pipher condition where one requires 
that $A$ be locally Lipschitz in $\Omega$, with
\begin{equation} \label{3a2}
\nabla A(X) \dist(X,E) \in L^{\infty}(\Omega),
\end{equation}
and 
that $|\nabla A(X)|^2 \delta(X)$ satisfy a Carleson measure condition, 
i.e.,  that there be a constant $C_M \geq 0$ such that
\begin{equation} \label{3a3}
\int_{\Omega \cap B(x,r)} |\nabla A(X)|^2 \dist(X,E) dX \leq C_M r^{n-1}
\end{equation}
for $x\in E$ and $0 < r < +\infty$.
As we shall see, we can manage with the weaker condition \eqref{3b3} because 
we don't need very precise estimates.

\smallskip
Our preferred geometric condition on $E$, uniform rectifiability, can also be defined using Definition~\ref{d1}, via the property called Bilateral Weak Geometric Lemma (BWGL).
The BWGL is known to imply apparently stronger rectifiability properties, such as the existence of
Big Pieces of Lipschitz Images of $\R^d$ inside of $E$; 
we refer to \cite{DS} for lots of information about uniform rectifiability, which is not the issue here. To define the WGL, we will use the following local version of the Hausdorff distance between two sets 
$E, F \subset \R^n$: we set 
\begin{multline}\label{1a8}
d_{x,r}(E,F) = 
\frac1r
\sup\big\{ \dist(y, F)\, ; \, y\in E \cap B(x,r) \big\} 
\\+ 
\frac1r
\sup\big\{ \dist(y, E)\, ; \, y\in F \cap B(x,r) \big\}
\quad\,
\end{multline}
for $x\in \R^n$ and $r > 0$; when $F \cap B(x,r) = \emptyset$, we take
$\sup\big\{ \dist(y, E)\, ; \, y\in F \cap B(x,r) \big\} = 0$, and similarly decide that
$\sup\big\{ \dist(y, F)\, ; \, y\in E \cap B(x,r) \big\} = 0$ when $E \cap B(x,r) = \emptyset$,
but we will probably not need to use this convention.

\begin{definition} \label{d2}
Let $E$ be an Ahlfors regular set of dimension $d$ in $\R^n$, where $d \in (0,n)$ is an integer.
We say that $E$ is uniformly rectifiable when for every $\varepsilon > 0$, the set
$\cG_{ur}(\varepsilon)$ of pairs $(x,r) \in E \times (0,+\infty)$ such that 
$d_{x,r}(E,P) \leq \varepsilon r$ 
for some $d$-plane $P = P(x,r)$, is Carleson-prevalent.
\end{definition}

As always with this type of conditions, we allow the Carleson constant for 
$\cB_{ur}(\varepsilon)$ to depend on $\varepsilon$ in any brutal way, 
but in the case of Definition \ref{d2}  
it turns out that we only need a single $\varepsilon$, chosen small enough,
depending on $n$, $d$, and the AR (Ahlfors regularity) constant for $E$.

\ms
Finally, turning to the properties of the Green function,
we shall find it more convenient to state our results in terms of $G^\infty$, 
the Green function with a pole at infinity. 
We will review the construction in Section \ref{S2} (see also \cite{DEM, DFM5}). 
Since our approximation property appears to be new, to motivate the forthcoming definition, 
let us first describe a situation which we judge perfect, and our conditions will try to measure 
how far we are from this situation.

In the case of $\Delta$ (or a constant coefficient operator), the perfect situation 
is when $d=n-1$ and $E$ is a $d$-plane. In this case the harmonic measure 
$\omega^\infty$ (with a pole at $\infty$) is a multiple of the Lebesgue measure on $E$, 
and the Green function $G^\infty$ is a multiple of the distance to $E$. 
In the case of our favorite operator $L_\alpha = - \dv D_\alpha^{d+1-n} \nabla$ 
from \eqref{1a2}, the perfect situation is when $d$ is an integer, $E$  is a $d$-plane, 
and $\mu$ is a multiple of the Hausdorff measure $\H^d_{\vert E}$ on $E$. 
In this case again the harmonic measure $\omega_L^\infty$ is a multiple of $\H^d_{\vert E}$ 
and the Green function $G^\infty$ is a multiple of the distance to $E$, which is also 
proportional to each $D_\beta$, $\beta > 0$. 

\begin{definition} \label{d3}
We say that $G^\infty$ is prevalently close to the distance to a $d$-plane when for each choice of
$\varepsilon >0$ and $M \geq 1$, the set $\cG_{Gd}(\varepsilon, M)$ of pairs 
$(x,r) \in E \times (0,+\infty)$ such that there exists a $d$-plane $P(x,r)$ and a positive constant 
$c > 0$, with
\begin{equation} \label{1a10}
|\dist(X,P) - c G^\infty(X)| \leq \varepsilon r
\text{ for } X \in \Omega \cap B(x,Mr),
\end{equation}
is Carleson-prevalent.

Similarly, given $\beta > 0$, we say that $G^\infty$ is prevalently close to $D_\beta$ when for each choice of
$\varepsilon >0$ and $M \geq 1$, the set $\cG_{GD_\beta}(\varepsilon, M)$ of pairs
$(x,r) \in E \times (0,+\infty)$ such that there exists a positive constant $c > 0$, with
\begin{equation} \label{1a12}
|D_\beta - c G^\infty(X)| \leq \varepsilon r
\text{ for } X \in \Omega \cap B(x,Mr),
\end{equation}
is is Carleson-prevalent.
\end{definition}

\smallskip
Here we find it convenient to have a constant $c > 0$ that we don't need to compute, 
especially since our function $G^\infty$ is only defined modulo a multiplicative constant. 
The definitions are easier to understand when we allow $M$ to be large, but in terms of 
Carleson-prevalent sets, taking $M=1$ and compensating with $\varepsilon$ would in fact 
give an equivalent result. The main positive result of this paper is as follows. 

\begin{theorem}\label{ti1} Assume that $\Omega\subset \RR^n$ is a domain with an unbounded $d$-Ahlfors regular boundary $E=\po$, $d$ integer. When $d=n-1$, assume, in addition, that $\Omega$ is 1-sided NTA. Let $L=-\dv A \nabla$ be an elliptic operator with bounded measurable coefficients when $d=n-1$, and more generally, given by 
$$L=-\dv A D_{\alpha, \mu}^{d+1-n} \nabla, \quad 0<d<n, $$
where $D_{\alpha, \mu}$ is the smooth 
distance \eqref{1a3} associated to some $\alpha>0$ and some AR measure $\mu$ on $E$, 
and $A$ is any elliptic matrix. 

If $G^\infty$ is prevalently close to the distance to a $d$-plane, in the sense of Definition~\ref{d3}, then $E$ is uniformly rectifiable. 

Let us assume, in addition, that $A$ is locally close to a constant coefficient operator when $d=n-1$, that is, \eqref{3b1}--\eqref{3b3} holds, or, when $d<n-1$, 
that $A$ is locally close to the identity matrix, 
that is, \eqref{3b1}--\eqref{3b3} holds with $A_0=I$. 

Then, conversely, $G^\infty$ is prevalently close to the distance to a $d$-plane whenever $E$ is uniformly rectifiable.
\end{theorem}

The theorem 
is a combination of Theorem~\ref{t3a1}, \ref{t3a1bis}, and \ref{t61}. As we discussed above, in both directions the conditions on the coefficients are of the nature of the best possible. The reader could perhaps be surprised that in domains with lower dimensional boundaries we require that $A$ is locally close to identity, rather than just a constant coefficient operator. This has to do with the underlying ``perfect" situation. Indeed, when $d=n-1$, the distance to a flat boundary is an affine function, which furnishes a solution for any constant coefficient operator. When $d<n-1$, there is a more delicate cancellation at place. Even in $\RR^n\setminus \RR^d$ 
the norm of $t$ (the component of $X$ in $\R^{n-d}$),
is a solution for $-\dv |t|^{d-n+1}\nabla$, but not necessarily for operators with a more general, even if constant, matrix of coefficients. Similarly, even if $E$ is a hyperplane, the distance to the boundary is only a solution for the emblematic operator \eqref{1a2}. One could maybe reach some more general results further generalizing the concept of a distance, but we chose not to pursue this direction in the present paper. 

Notice that when we prove that $G$ is prevalently close to $\dist(X,P)$, for instance,
this does not imply that it is Lipschitz, because the constant $c$ in \eqref{1a10} may
depend on $B(x,r)$.
For instance, if $E$ is a Lipschitz domain with a corner at $0$, we know that $G$ may behave 
like $\dist(X,E)^\alpha$, $\alpha \neq 1$, near $0$. The balls centered at the origin
probably lie in the complement of the good set $\cG_{Gd}(\varepsilon, M)$, 
but there are many good balls $B(y,t)$ near $0$, where $E$ is flat and \eqref{1a10} 
or \eqref{1a12} holds, and for these balls the constants $c= c(y,t)$ will typically tend to 
$0$ or $+\infty$ as $(y,t)$ tends to $(0,0)$.

Going further, let us discuss the ``weak" nature of the results. 
The way we formulate our theorems, via a weak 
Carleson packing condition, or one could also say, without a precise control of the constants, 
is not necessarily optimal, or rather, looks stronger than expected in one direction and weaker than expected in the other. 
This has to do with a self-improvement of scale invariant estimates which takes place in 
uniformly rectifiable sets. Indeed, in the context of uniform rectifiability for AR sets,
the Bilateral Weak Geometric Lemma of Definition \ref{d2} is known to be equivalent 
to stronger and more precise definitions, for instance, the aforementioned existence of 
big pieces of Lipschitz images uniformly at all scales \cite{DS}. 
Similarly, here when we say that $G^\infty$ is prevalently close to $\dist(X,P)$, 
we prove that
$\cG_{GD}(\varepsilon,M)$ is prevalent for each choice of $M$, $\varepsilon$, but we do not give the rate of convergence of the Carleson packing constant for
$E \times (0,+\infty) \sm \cG_{GD}(\varepsilon,M)$ in terms of $\varepsilon$.
Ultimately, we expect square Carleson estimates on the quantities that control $\dist(X,P) - c G^\infty$,
and intend to prove such estimates in incoming papers with L.~Li and J.~Feneuil. In fact, the strong estimates would pertain to the ($L^2$ averages of) the gradient of the Green function, rather than the Green function itself, and would necessarily have to be properly localized. With this in mind, we are treating the derivatives of $G^\infty$ as well and then proving the corresponding local estimates in Sections~\ref{S4} and \ref{S5}, respectively.

\smallskip

For now we have only discussed approximation of the Green function by 
the distance to a plane,
that is, \eqref{1a10}, and only in integer dimensions. The second main question in the present paper concerns the  approximation by a distance to the initial boundary $E$, as in \eqref{1a12}. Despite an apparent similarity to \eqref{1a10}, it is actually more intricate at least in one direction, for even the fact that $G^\infty$ is on the spot equal to a multiple of the distance to $E$, $D_\beta$, might or might not ensure that the boundary is nice. Even the fact that the dimension is integer has to be proved and, in some circumstances, can fail. For brevity, we will formulate all results in the context of operators which are locally close to the constant coefficient ones, even though some results are still true in the full generality of all elliptic operators -- see Section~\ref{S7}. We start with the classical case. 

\begin{theorem}\label{ti2} Assume that $\Omega\subset \RR^n$ is a domain with an unbounded $d$-Ahlfors regular boundary $E=\po$, $n-2<d<n$. When $d\geq n-1$, assume, in addition, that $\Omega$ is 1-sided NTA. Let $L=-\dv A \nabla$ be an elliptic operator with bounded measurable coefficients which are locally close to constant coefficient matrices in the sense of \eqref{3b1}--\eqref{3b3}.

If $E$ is uniformly rectifiable, then for every choice of $\beta > 0$ and any AR measure
$\mu$ on $E$, $G$ is prevalently close to $D_{\beta,\mu}$.

Conversely, if the Green function $G^\infty$
is prevalently close to $D_{\alpha,\mu}$ for some $\alpha > 0$ and some AR  
measure $\mu$ of dimension $d$ on $E$ then $d$ is integer and  $E$ is uniformly rectifiable.
\end{theorem}

This is a combination of Theorems~\ref{t3a1}, \ref{t71}, and Proposition~\ref{t72}, 
the latter being valid for all elliptic operators. 

The situation for domains with lower dimensional boundaries is trickier, due to the aforementioned rather mysterious fact that for a certain very special choice of coefficients the situation could be ``perfect" in terms of PDEs without any regularity or flatness. 
This is one of the main discoveries in \cite{DEM}, and since that paper, 
we refer to it as a ``magic $\alpha$" case. 
Indeed, if $E$ is any unbounded $d$-Ahlfors regular set and $\alpha=n-d-2$, 
$R_{\alpha}$ is just the convolution of $\mu$ with a multiple of the Laplace's fundamental 
solution in $\R^n$, so it is harmonic, and a direct computation shows that 
$L_\alpha D_{\alpha, \mu}=0$ in $\Omega$ for $L_\alpha$ given by \eqref{1a2}. 
Hence, by uniqueness, $D_{\alpha, \mu}$ is a multiple of $G^\infty$, and $C^{-1}\mu \leq \omega^\infty \leq C \mu$. In fact, 
if $d$ is an integer, and $E$ is rectifiable (not necessarily uniformly rectifiable)
with $\mu = \H^d_{\vert E}$, then the harmonic measure $\omega^\infty$ with pole at 
infinity is proportional to $\mu$. In all the other cases, we only have that
$C^{-1}\mu \leq \omega^\infty \leq C \mu$, essentially because we cannot say that
$\frac{\d \omega^\infty}{\d \mu}$ is the normal derivative of $G^\infty$ when $E$ is unrectibiable, but we can argue that getting the measure $\mu$ point blank makes less sense
in this case, because the density of $\mu$ does not exist.
Thus, for such a ``magic $\alpha$" 
the  
elliptic measure is absolutely continuous with respect to the Hausdorff measure, with a  
density given by an $A^\infty$ weight, on any Ahlfors regular set.
Notice that we need $d < n-2$ for this to happen as $\alpha>0$, but we do not need $d$ to be an integer, and certainly we do not need $E$ to be uniformly rectifiable. This is a strangely degenerate case which does not resonate with anything we know about the standard domains with $n-1$ dimensional boundaries, but it also turns out to be immensely useful as one basically gets to use  the distance $D_\alpha$ as the Green function, getting the best of both worlds: a solution to a PDE and an explicit, easy-to handle formula. This is a rare luxury as typically we do not know explicitly the Green function -- see, e.g., \cite{Fen} for some applications. However, in the context of this paper this case is certainly not amenable to any free boundary results akin to Theorem~\ref{ti2}. We hope, however, that this is an isolated miraculous cancellation and generally  the fact that $G^\infty =C D_{\alpha, \mu}$  implies that the set is flat and $\mu$ is a multiple of a flat measure. Not to extend the introduction any further, let us refer the reader too Section~\ref{S8} for a detailed statement of this condition, which we will denote $\Upsilon_{\rm{flat}}$.

\begin{theorem}\label{ti3} Assume that $\Omega\subset \RR^n$ is a domain with an unbounded $d$-Ahlfors regular boundary $E=\po$. When $d\geq n-1$, assume, in addition, 
that $\Omega$ is 1-sided NTA. Let 
$$L=-\dv A D_{\alpha, \mu}^{d+1-n} \nabla, \quad 0<d<n, $$
where $D_{\alpha, \mu}$ is the smooth distance 
distance \eqref{1a3} associated to some $\alpha>0$ and some AR measure $\mu$ on $E$, 
and $A$ is any elliptic matrix with bounded measurable coefficients which are locally close to the identity in the sense that \eqref{3b1}--\eqref{3b3} holds with $A_0=I$.

If $E$ is uniformly rectifiable, then for every choice of $\beta > 0$ and any AR measure
$\nu$ on $E$ (possibly different from $\mu$), $G$ is prevalently close to $D_{\beta,\nu}$.

Conversely, if the condition 
$\Upsilon_{\rm{flat}}(d,\alpha,\Delta)$ of Definition \ref{d71} 
holds and the Green function $G^\infty$
is prevalently close to $D_{\alpha,\mu}$ then $d$ is integer and $E$ is uniformly rectifiable.

\end{theorem}

This is a combination of Theorems~\ref{t3a1bis} and Theorem~\ref{t81}. One can find more related results and an extended discussion of the condition $\Upsilon_{\rm{flat}}$ in Section~\ref{S8}. For now, let us wrap up the introduction and send an interested reader to the body of the paper.

\section{The Green function at infinity and the basic result about limits}
\label{S2}
We first remind the reader of how the Green function $G^\infty$ is constructed.
We have two main cases in mind, which we rapidly describe now.  Let us start with the geometric assumptions.

Throughout the paper, $\Omega$ is a domain in $\RR^n$ such that $E = \d\Omega$ is an unbounded Ahlfors regular set of dimension $d \in (0,n)$, so that \eqref{1a1} is satisfied. There are natural counterparts of our results for bounded domains, but for simplicity of notation we will concentrate on the unbounded case.  

When $d \geq n-1$, we also demand that $\Omega$ is a uniform ({\it{aka}} a one-sided NTA) domain, i.e., that
\begin{equation} \label{2b4}
\text{$\Omega$ has interior corkscrew points and Harnack chains.} 
\end{equation}
Let us say what this means.
First, we require the existence of corkscrew points: for $x\in E$ and $r>0$, 
we demand the existence of a point $A = A_{x,r}$ such that
\begin{equation} \label{2a1}
A \in \Omega \cap B(x,r) \ \text{ and } \dist(A,E) \geq C^{-1} r
\end{equation}
(we shall call $A$ a corkscrew point for $B(x,r)$).

Secondly, we require 
that for each $M \geq 1$, we can find an integer $N =N(M)$ such that when 
$X, Y \in \Omega$ are such that $|X-Y| \leq M \min(\dist(X,E),\dist(Y,E))$, 
we can find a (Harnack) chain of balls $B_1,  \ldots, B_N$ 
such that $B_1$ contains $X$, $B_N$ contains $Y$, each $B_i$, 
$2 \leq i \leq N$, meets $B_{i-1}$, and the $2B_i$ are all contained in $\Omega$. 
Let us not play useless games here;  
we also require that all the $B_i$ of the chain have a 
radius at least $C^{-1} \min(\dist(X,E),\dist(Y,E))$.

Another way to formulate the Harnack chain condition would be to require the existence of a path
in $\Omega$ that goes from $X$ to $Y$, stays at distance at least 
$C^{-1} \min(\dist(X,E),\dist(Y,E))$ from $E$, and has a length at most $C|X-Y|$.

We shall refer to $C$ in \eqref{1a1} as an AR constant of $E$, and to $C$ and $M$ is the corkscrew and Harnack chain conditions above as 1-sided NTA constants.

A uniform domain with an Ahlfors regular boundary is sometimes referred to as a 1-sided chord-arc domain, but we will rarely use this terminology. 

The interior corkscrew and Harnack chain conditions are merely quantitative openness and path connectedness of the domain $\Omega$. They guarantee a 
reasonable behavior of the Green function at the level of fundamental estimates and so we prefer to assume that the domains are uniform when $d\geq n-1$ throughout the paper. When $d < n-1$, it turns out that interior corkscrew and Harnack chain conditions are automatically satisfied for any domain with an Ahlfors regular boundary $E=\po$ of dimension $d$, essentially because $E$ is so thin. This is checked in \cite{DFM2}. 

\ms

Let us turn to 
the elliptic operators considered in the present paper. We split into two cases, essentially corresponding to the dimension of the boundary of the set, although they are not completely exclusive. 

In {\bf Case 1}, which we will also call the classical case, we are given a domain $\Omega$ with $E=\po$, an unbounded Ahlfors regular set $E$ of dimension $d$, with $n-2 < d < n$. 
When $d\geq n-1$, we assume, in addition, that the domain is uniform, i.e., the interior corkscrew and Harnack chain conditions are satisfied. 
We consider an elliptic divergence form operator
\begin{equation} \label{2b1}
L = - \dv A \nabla,
\end{equation}
where the matrix $A$ of coefficients is measurable and satisfies the usual ellipticity properties
\begin{equation} \label{2b2}
|\langle A(X) \xi, \zeta \rangle| \leq C_e |\xi||\eta| \, \text{ for $X\in \Omega$, $\xi \in \R^n$, and $\zeta \in \R^n$,} 
\end{equation}
\begin{equation} \label{2b3}
\langle A(X) \xi, \xi \rangle \geq  C_e^{-1} |\xi|^2 
\, \text{ for $X\in \Omega$ and $\xi \in \R^n$,} 
\end{equation}
and some constant $C_e \geq 1$.
For some results we will also require $L$ to be close enough to constant 
coefficient operators, but we will explain this in due time.

In {\bf Case 2}, which is not exclusive of Case 1, $\Omega$ is still a domain in $\RR^n$
such that $E = \d\Omega$ is an unbounded Ahlfors regular of dimension $d$, but now 
we allow any dimension $d \in (0,n)$. When $d \geq n-1$, we require that $\Omega$ is uniform, i.e., that  \eqref{2b4} is satisfied, as before. 
And now we consider our favorite operator
\begin{equation} \label{2b6}
L = L_\alpha = - \dv D_{\alpha,\mu}^{d+1-n} \nabla
\end{equation}
of \eqref{1a1}, where $\alpha > 0$ and $D_{\alpha,\mu}$ is given by \eqref{1a3} and \eqref{1a4}, for some AR measure $\mu$ on $E$, or its generalization 
\begin{equation} \label{2b6-bis}
L =  -\dv A\,D_{\alpha,\mu}^{d+1-n} \nabla,
\end{equation}
where $A$ satisfies \eqref{2b2}--\eqref{2b3}. Clearly, at such a level of generality we could as well say $L =  -\dv A\,\dist^{d+1-n} 
(\cdot, E)\nabla,$ but writing $L$ as above will occasionally be somewhat more convenient. Clearly, so far Case 1 and Case 2 coincide when $d=n-1$ and in fact, we could further generalize Case 2 to cover Case 1 completely, but we prefer to keep them separate for reasons that will become evident a little later: the additional assumptions down the road will start to deviate.

\ms

Let us review the definition of the Green function with a pole at $\infty$,
and then state the main result about limits that will be behind all our proofs by compactness.

First recall that in both cases, we can associate to $\Omega$ and $L$ a collection of probability measures
$\omega^X = \omega^X_L$, referred to as the harmonic (or elliptic) measure with pole at $X \in \Omega$,
and that have standard properties (such as doubling) that will be recalled when we need them.
There are quite a few papers listing these fundamental estimates in Case 1, at least when $d=n-1$ (and perhaps more generally). Essentially, one can say that the program outlined in \cite{JK} for 2-sided NTA domains still applies. In Case 2, this is done in \cite{DFM2} when $d < n-1$. 
Finally, both cases for all relevant dimensions are covered in \cite{DFM5}.
See Section 3 of \cite{DFM5} where it is explained how the present assumptions 
are covered by the 
hypotheses of \cite{DFM5}. The setting of \cite{DFM5} is much more general than what we need here, but it is convenient to have all the relevant references in one place and so we will mainly use \cite{DFM5} for fundamental elliptic theory: Caccioppoli, De Giorgi-Nash-Moser, existence and uniqueness of elliptic measure, Green functions, etc. 
In particular, the reader may find in the same source 
a construction of the Green function $G^X$ for $L$, 
at the pole $X\in \Omega$. Again these functions enjoy the usual properties of Green functions
in NTA domains, and we will recall these properties (and those of solutions of $Lu=0$) when we need them.

\ms
Next we say a few words about the Green function $G^\infty$ with pole at $\infty$, which we construct as a limit of functions $G^X$. 
We will also use the opportunity to recall some notation and estimates that will be used later.

First of all, we define a weight $w$ on $\Omega$ by
\begin{equation} \label{2b9}
w(X) = 1 \text{ in Case 1 and $w(X) = \dist(X,E)^{d+1-n}$ in Case 2.}
\end{equation}
Then we let $W=W(\Omega)$ denote the Hilbert space of functions $u \in L^1_{loc}(\Omega)$ whose derivative
(in the sense of distributions) lies in $L^2(\Omega, w(X)dX)$; thus
\begin{equation} \label{2c8}
||u||_W^2 = \int_\Omega  |\nabla u(X)|^2 w(X) dX
\end{equation}
is finite.  We also need local versions $W_r(B)$ of $W$, which are defined in Section 8 of \cite{DFM2} or Section 10 of \cite{DFM5}, as follows. Given an open set $H'\subset \RR^n$ and $H:=H'\cap \overline\Omega$, we let 
$$
W_r(H) := \{u \in L^1_{loc}(H \cap \Omega,m):\, \varphi u \in W \text{ for all } \varphi \in C^\infty_0(H')\}.
$$
It is natural to call this space $W_r(H)$, as opposed to $W_r(H')$, because it does not depend 
on the part of $H'$ that lives away from $\overline\Omega$. It does or does not carry information about the behavior of $u$ near the boundary depending on whether $H$ is properly contained in the domain $\Omega$. Clearly, $\nabla u \in L^2(B, w(X) dX)$ when 
$u \in W_r(B)$. The Green function $G^X$ with pole at $X \in \Omega$ lies in $W_r(B)$ for any ball
$B\subset \RR^n$ such that $X \notin 2B$ (see Lemma 10.2 in \cite{DFM2} or Lemma 14.60 in \cite{DFM5}).

We say that the function $u$ is  a weak solution to $Lu=-\div A \nabla u=0$ in $\Omega$ when it lies in all the local spaces $W_r(B)$ and 
\begin{equation} \label{2a7}
\int_{\Omega} A \nabla u \cdot \nabla \varphi = 0
\end{equation}
for every $\varphi \in C^\infty_0(\Omega)$. Here, $A$ can be as in \eqref{2b2}--\eqref{2b3} in Case 1 or \eqref{2b6}--\eqref{2b6-bis} in Case 2. 

\ms

We are now ready to describe how we construct the Green function $G^\infty$ with pole at $\infty$.
Fix any ball $B_0 = B(x_0, r_0)$ centered on $E$,
choose a corkscrew point $A_0$ for $B_0$
(see \eqref{2a1}), take any sequence $\{ X_k \}$ in $\Omega \sm B_0$ such that 
$\lim_{k \to +\infty} |X_k| = +\infty$, and consider the functions
\begin{equation} \label{2b7}
g_k(Y) = \frac{G^{X_k}(Y)}{G^{X_k}(A_0)}.
\end{equation}
Notice that $g_k$ is nonnegative, $L$-harmonic on $\Omega \cap B(x_0, |X_k-x_0|/2)$, 
and $g_k \in  W_r(B)$ for every ball $B$ with 
constants that are uniform in $k$,
as soon as $k$ is so large that $X_k \notin 2B$.
It also has a vanishing trace at the boundary, by definition of any $G^X$.

By Harnack's principle and the normalization $g_k(A_0) = 1$, we see that for each compact 
set $K \subset \Omega$, there is a constant $C_K \geq 1$ such that
\begin{equation} \label{2b7a}
C_K^{-1} \leq g_k(Y) \leq C_K \text{ for $Y \in K$},
\end{equation}
as soon as $k$ is large enough (depending on $K$), so that $X_k$ no longer lies
in the union of a (finite) collection of Harnack chains that connect $A_0$ to any point of $K$.

Similar estimates hold near $E$, i.e., on $\Omega \cap B$, where $B$ is any ball centered
on $x_0\in E$. Indeed, the functions $g_k$ are continuous, and even H\"older continuous, at the boundary. With this in mind, let us extend $g_k$ by zero to the rest of $\RR^n$.
We use Lemma 11.50 in \cite{DFM2}
or Lemma 15.14 in \cite{DFM5}, which says that $g_k(Y) \leq C g_k(A^B)$ for
$Y \in B$ (and $k$ large enough), where $A^B$ is a corkscrew point for $B$.
But $g_k(A^B)$ is controlled by \eqref{2b7a}, so we get that for each large ball $B$ centered
at $x_0$, there is a constant $C_B$ such that
for $k$ large,
\begin{equation} \label{2c11}
g_k(Y) \leq C_B \text{ for $Y \in B$}.
\end{equation}
Then we can use the Caccioppoli estimate at the boundary (Lemma 8.47 in \cite{DFM2}
or Lemma 11.15 in \cite{DFM5}), and get that for $B$ as above (and $k$ large)
\begin{equation} \label{2c12}
\int_{B\cap \Omega} |\nabla g_k|^2 w(X)dX 
\leq C(B) \int_{ C B\cap \Omega} |g_k|^2 w(X)dX \leq C(B).
\end{equation}
Here, as usual the fact that $g_k$ vanishes on $E$ and $g_k \in W_r(B)$ for any given $B$ 
(and for $k$ large) was used to check the assumptions.
Finally, we can also use the H\"older continuity at the boundary (Lemma 8.106 in \cite{DFM2} 
or Theorem 1.6 in \cite{DFM5}) to show that 
\begin{equation} \label{2c13}
\text{$g_k$ is H\"older continuous on $B$, with exponent $\beta$ and constant $C(B)$.}
\end{equation}
Of course the important point for the moment is that none of the constants $\beta, C_B,$ and $C(B)$ depends on $k$, even though they are only valid for $k$ large, depending on $B$.

Return to the $g_k$. Because of \eqref{2c13} it is easy to extract a subsequence so that 
$\{ g_k \}$ converges uniformly on compact subsets of $\RR^n$ to a limit $G^\infty$. 
In fact, using again the Caccioppoli estimate at the boundary and passing to a subsequence, 
we also get that for each ball $B$,  $\nabla g_k$ converges weakly in $L^2(\Omega \cap B, w)$ 
to $\nabla G^\infty$, and then that $G^\infty$ is a (weak) solution of $L$ (we will showcase the details in similar, but more delicate, arguments soon). Since $g_k$ is  
H\"older continuous and vanishing on the boundary, 
$G^\infty$ also vanishes at the boundary. And finally $G^\infty \in W_r(B)$ for every $B$ 
(using the definition of $W_r$, weak convergence, and the estimates above).

Now, $G^\infty$ is also unique, modulo a multiplicative constant. That is, if $G$ is positive 
on $\Omega$, lies in all the spaces $W_r(B)$, it is $L$-harmonic on $\Omega$, and has a 
vanishing trace on $E$ (the latter is well defined when $G\in W_r(B)$), then $G = c G^\infty$ 
for some $c > 0$. This follows from the comparison principle 
(Theorem 11.146 in \cite{DFM2} or Theorem 1.16 in \cite{DFM5}), 
plus some algebraic manipulations on the oscillation of $G/G^\infty$ (similar to what one does for
the H\"older continuity at the boundary) that the reader may find in \cite{DEM}, 
Corollary 6.4 and Lemma~6.5. This procedure coherently and uniquely defines $G^\infty$
(modulo a multiplicative constant),
which we will refer to as the Green function with the pole at infinity.

\ms
We are ready for our main result about limits. 
We assume that we have a sequence of open sets $\Omega_k$, 
bounded by Ahlfors regular sets $E_k$, and operators $L_k$, that all satisfy the assumptions 
of this paper (for either Case 1 or Case 2), with uniform estimates. 

We assume (for convenience) that all the $E_k$ contain the origin, and that 
$\{ E_k \}$ converges to a closed set $E_\infty$ and 
$\{ \Omega_k \}$ converges to an open set $\Omega_\infty$, in the sense that with the notation of 
\eqref{1a8}, 
\begin{equation} \label{2a2}
\lim_{k \to +\infty} d_{x,r}(E_\infty, E_k) = 0
\ \text{ and } \ 
\lim_{k \to +\infty} d_{x,r}(\Omega_\infty, \Omega_k) = 0
\end{equation}
for every choice of $x\in \R^n$ and $r >0$ (we may also restrict to $x\in E_\infty$;
this would be equivalent).

We shall check soon that $E_\infty$ is the boundary of $\Omega_\infty$, 
$E_\infty$ is Ahlfors regular of dimension $d$, and $\Omega_\infty$ satisfies our 
one-sided NTA condition \eqref{2b4}, but let us continue our description.

We also need the operators $L_k = - \dv A_k \nabla$ to converge to a limit 
$L_\infty = -\dv A_\infty \nabla$ (and similarly for the operators of Case 2), and we require that
\begin{equation} \label{2a10bis}
\lim_{k \to +\infty} ||A_k - A_\infty||_{L^1(B)} = 0 
\ \text{ for every ball $B$ such that $2B \subset \Omega_\infty$.}
\end{equation}
Our proof would be a little simpler under the stronger local $L^\infty$ convergence 
of coefficients, namely when
\begin{equation} \label{2a10}
\lim_{k \to +\infty} ||A_k - A_\infty||_{L^\infty(B)} = 0 \ \text{ for every ball $B$ such that $2B \subset \Omega_\infty$,}
\end{equation}
which would probably be more reasonable if we wanted better quantitative results, 
but the fact that we are only interested in weak results allows us to use the 
weaker condition \eqref{2a10bis}.

Finally we choose a corkscrew point $A_0$ (relative to $\Omega_\infty$) for some ball 
$B_0$ centered on $E_\infty$, which will be used to normalize the Green functions.

 The next Theorem shows that with the assumptions above, the functions $G_k^\infty$ 
 associated to the $\Omega_k$ and the $L_k$ and normalized by $G_k^\infty(A_0) = 1$ 
 converge to the Green functions $G_\infty^\infty$ associated to $\Omega_\infty$ and $L_\infty$ 
 and normalized by $G_\infty^\infty(A_0) = 1$.

\begin{theorem} \label{l2a4}
Let the $\Omega_k$ be
domains in $\RR^n$ with unbounded $d$-dimensional Ahlfors regular boundaries $E_k=\po_k$, $d\in (0, n)$, corresponding to AR measures $\mu_k$, and, in addition, satisfying the interior Harnack chain and corkscrew conditions \eqref{2b4} when $d\geq n-1$, with all AR and 1-sided NTA constants uniform in $k$. 
Assume that the origin belongs to all the $E_k$,
that the domains $E_k$ converge to a closed set $E_\infty$, and that the 
$\Omega_k$ converge to an open set $\Omega_\infty$ in the sense of \eqref{1a2}. 
Then $E_\infty$ is the boundary of $\Omega_\infty$, $E_\infty$ is Ahlfors regular of dimension $d$,
and $\Omega_\infty$ is uniform, i.e., satisfies \eqref{2b4}.
If $\mu_k$ is an AR measure on $E_k$ (with uniform bounds) and $\mu$
is any weak-* limit of the $\mu_k$, then $\mu$ is an AR measure on $E_\infty$.

Assume furthermore that either 
\begin{equation} \label{218c1}
L_k = - \dv A_k \nabla, \quad L_\infty = - \dv A_\infty \nabla, \quad n-2<d<n, 
\end{equation}
or 
\begin{equation} \label{218c2}
L_k = - \dv A_k D_{\alpha, \mu_k}^{d+1-n} \nabla, \quad L_\infty 
= - \dv A_\infty D_{\alpha, \mu}^{d+1-n} \nabla, \quad 0<d<n
\end{equation}
(with $\mu_k \rightharpoonup \mu$ as above),
in both both cases subject to \eqref{2a10bis}.

Then the functions $G_k^\infty$ associated to the $\Omega_k$ and the 
$L_k$ and normalized by $G_k^\infty(A_0) = 1$ converge, uniformly on every compact 
subset of $\RR^n$, and in $W^{1,2}_{loc}(\Omega_\infty)$, 
to the Green functions $G_\infty^\infty$ associated to $\Omega_\infty$ and $L_\infty$ and normalized by $G_\infty^\infty(A_0) = 1$.
\end{theorem}

Here and everywhere, we say that the sequence of Radon measures $\mu_k$ on $\RR^n$ weak-* converges to some Radon measure $\mu$, and write $\mu_k \rightharpoonup \mu$, if 
$$ \int f\,d\mu_k \to \int f \,d\mu \quad \mbox{for any} \quad f\in C_c(\RR^n).$$ 
Notice that there always exist weak-* limits $\mu$, as in the statement.

As usual, \eqref{218c1} and \eqref{218c2} will be referred to as Case 1 and Case 2, respectively.

\begin{remark} \label{r2a18-bis}
Formally speaking, the Green functions are only defined in their corresponding domains, but using their H\"older continuity inside and at the boundary, we silently extend them by zero to the entire $\RR^n$, and hence we can justifiably talk about the uniform convergence on compacta in $\RR^n$.
\end{remark}

\begin{remark} \label{r2a18}
The same statement also holds if we replace the Green functions $G_k^\infty$
with Green functions $G_k = c_k G_k^{Y_k}$ computed at poles $Y_k$ and normalized
so that $G_k(A_0) = 1$, if we also assume that $\lim_{k \to +\infty} |Y_k| = + \infty$.
We will track this case too along the proof.
\end{remark}

\noindent {\it{Proof of Theorem~\ref{l2a4}}}. Let us first check that $(\Omega_\infty, E_\infty)$ satisfies our desired geometric properties.
We shall only highlight some elements of the proof, as the reader can consult \cite{HMMTZ1} for the details in the case of $d=n-1$ which transfer virtually verbatim to our more general  situation. 

The fact that $\d\Omega_\infty = E_\infty$ will be easy. 
Maybe we should notice first that \eqref{2a2} also implies that 
$d_{x,r}(\R^n \sm \ol \Omega_k, \R^n \sm \ol \Omega_\infty)$
also tends to $0$ for all pairs $(x,r)$: for instance, if $z \in \R^n \sm \ol \Omega_k$, 
it cannot be far from $\R^n \sm \ol \Omega_\infty$ because otherwise it lies in the middle 
of $\Omega_\infty$, hence also of $\Omega_k$). 
Next, if $z \in \d\Omega_\infty$, then for each $\varepsilon >0$ the ball
$B(z,\varepsilon)$ meets both $\Omega_\infty$ and $\R^n \sm \Omega_\infty$.
Then for $k$ large, $B(z,\varepsilon)$ also meets $\Omega_k$ and $\R^n \sm \Omega_k$,
hence also $E_k$. Now $\{ E_k \}$ converges to $E_\infty$, and $E_\infty$ is closed; 
it follows that $z\in E_\infty$. 
Finally, let $z\in E_\infty$ be given; for each $\varepsilon > 0$ we can find points 
$z_k \in E_k \cap B(z,\varepsilon)$ for $k$ large, hence also points 
$x_k \in \Omega_k \cap B(z,\varepsilon)$ and $y_k \in B(z,\varepsilon) \sm \ol\Omega_k$,
and by the extension of \eqref{2a2}, points $x'_k \in \Omega_\infty \cap B(z,2\varepsilon)$
and $y'_k \in B(z,2\varepsilon) \sm \ol\Omega_\infty$. 
So $z\in \d\Omega_\infty$ and $\d\Omega_\infty = E_\infty$. 

From here, we can show that $E_\infty$ supports a measure which is a weak limit of $\mu_k$ and which is Ahlfors regular too. Moreover, \eqref{2b4} also holds for $\Omega_\infty$.
For instance, if $X, Y \in \Omega_\infty$ are given, then for $k$ large we can find 
$X_k, Y_k \in \Omega_k$, as close as we want to $X$ and $Y$, and a Harnack chain 
for $X_k$ and $Y_k$ in $\Omega_k$ will also work, with very minor modifications, 
for $X$ and $Y$ in $\Omega_\infty$.
Again, the reader can consult \cite{HMMTZ1} for more details. For us, the main point of this verification was that we can apply the results 
of \cite{JK, DFM2, DFM5} to $L_\infty$ and $\Omega_\infty$. 
In particular, we can talk about $G_\infty^\infty$ and use the above uniqueness result for $G_\infty^\infty$.

We are ready for the PDE part of the argument. 
Set $G_k = G_k^\infty$ for all $k$ to save notation. We can run the same limiting argument  as the one which allowed us to define $G^\infty$ on a given domain. Indeed, \eqref{2c11} still holds in any ball $B\subset \RR^n$ (not necessarily contained in $\Omega_\infty$) with a constant uniform in $k$ as long as $G_k$ is extended by zero to the complement of $\Omega_k$.  (If we take a limit of  $G_k = c_k G_k^{Y_k}$  as in Remark~\ref{r2a18} instead, we also have to make sure that $k$ is large enough for $Y_k$ to stay away from $B$, but the constants are still uniform in $k$). Furthermore, the analogue of \eqref{2c12} holds: 
\begin{equation} \label{2c12b}
\int_{B\cap \Omega_k} |\nabla G_k|^2 w_k(X)dX 
\leq C(B) \int_{ C B\cap \Omega_k} |G_k|^2 w_k(X)dX \leq C(B),
\end{equation}
with uniform constants, which again yields \eqref{2c13}. Here, $w_k$ is as in \eqref{2b9}.
Thus, much as before, we can extract a subsequence which converges uniformly on compacta 
in $\RR^n$ to some $G_\infty$, continuous and equal to zero on the complement of 
$\Omega_\infty$, and such that $\nabla G_k$ converges weakly to $\nabla G_\infty$ in $L^2_{loc} (\Omega_\infty)$ (notice that the weights are irrelevant when we stay away from the boundary). Eventually, we will prove the strong convergence  of $G_k$ 
in $W^{1,2}_{loc} (\Omega_\infty)$, but for now let us continue.

The function $G_\infty$ is in $W_r(B)$ for any ball $B\subset \RR^n$. Indeed, it is sufficient to consider $B$ centered on $E_\infty$.
For a single $G_k$, \eqref{2c12b} holds with uniform in $k$ constants.
Then we can look at the contribution to \eqref{2c12b} of any fixed compact subset $H$ 
of $\Omega_\infty \cap B$, let $k$ tend to $+\infty$ in \eqref{2c12b}, and get that 
\begin{equation} \label{2c39}
\int_{H} |\nabla G_\infty|^2 w_\infty(X)dX \leq C(B),
\end{equation}
for instance by Fatou, or the weak convergence in $W^{1,2}_{loc}(\Omega_\infty)$. 
Then we use the fact that $C(B)$ does not depend on $H$, take a supremum, and get that
\begin{equation} \label{2c40}
\int_{\Omega_\infty \cap B} |\nabla G_\infty|^2 w_\infty(X) dX \leq C(B).
\end{equation}
This (together with the local H\"older continuity of $G$ and the local integrability 
of the weight $w_\infty$) is actually stronger than the fact that $G \in W_r(B)$ for every ball, which requires that $\varphi G$ lies in the Hilbert space $W$ of \eqref{2c8} 
for every $\varphi \in C_0^\infty$ (see the beginning of Section 8 in \cite{DFM2}).

We also need to show that $G$ is a weak solution of $L_\infty G = 0$.
In Case 1,  we write
\begin{equation} \label{2c41}
\int_{\Omega} A_\infty \nabla G_\infty \cdot \nabla \varphi = \int A_\infty (\nabla G_\infty-\nabla G_k) \cdot \nabla \varphi + \int (A_\infty-A_k) \nabla G_k \cdot \nabla \varphi + \int  A_k  \nabla G_k \cdot \nabla \varphi
\end{equation}
for every $\varphi \in C^\infty_0(\Omega_\infty)$ and any $k$ such that $\Omega_k$ contains the support of $\vp$.  The first integral converges 
to zero by the weak convergence of $\nabla G_k$ in $L^2_{loc}(\Omega_\infty)$, 
the second one is bounded by $\sup |\nabla \vp| \int_{\supp \vp} |A_\infty-A_k| |\nabla G_k|$, 
and hence, it converges to zero by \eqref{2a10bis} and \eqref{2c12b}, and the last one is zero because $G_k$ is a solution to $L_k$. Hence, $G_\infty$ is a solution to $L_\infty$, as desired. 

A similar argument also handles Case 2. Indeed, we can write \eqref{2c41} with $A_\infty D_{\alpha, \mu}^{d+1-n}$ in place of $A_\infty$ and  $A_k D_{\alpha, \mu_k}^{d+1-n}$ in place of $A_k$. The weights are harmless since we are away from the boundary, we only have to show that $A_\infty D_{\alpha, \mu_k}^{d+1-n}$ converge to $A_k D_{\alpha, \mu}^{d+1-n}$ to handle the analogue of the second term on the right-hand side of \eqref{2c41}. However, $D_{\alpha, \mu_k}^{d+1-n}$ converge to $D_{\alpha, \mu}^{d+1-n}$ uniformly on any compact set in $\Omega_\infty$. Indeed, on any set away from the boundary it is sufficient to show that $R_{\alpha, \mu_k}$ converge to $R_{\alpha, \mu}$ (see \eqref{1a4}) which follows from the weak-* convergence of $\mu_k$ 
to $\mu$ and the Ahlfors regularity of the measures $\mu_k$ and $\mu$ (the latter allowing us to restrict attention to a compactly supported approximation of $|X-y|^{-d-\alpha}$). Hence, writing 
\begin{equation}
\label{eqConv}
A_\infty D_{\alpha, \mu}^{d+1-n}- A_k D_{\alpha, \mu_k}^{d+1-n}=(A_\infty-A_k)D_{\alpha, \mu_k}^{d+1-n}+A_\infty (D_{\alpha, \mu}^{d+1-n}-D_{\alpha, \mu_k}^{d+1-n}),
\end{equation}
we see that \eqref{2a10bis} holds with $A_\infty D_{\alpha, \mu}^{d+1-n}$ in place of $A_\infty$ and  $A_k D_{\alpha, \mu_k}^{d+1-n}$ in place of $A_k$.

So $G$ is a solution for $L_\infty$ that vanishes on $E_\infty$, and it lies in the 
correct spaces $W_r(B)$. Since $\Omega_\infty$ and $A_\infty$ satisfy the same 
assumptions as the $\Omega_k$ and $A_k$, we can apply the result of uniqueness for the 
Green function, and we find out that $G$ is a Green function for $L_\infty$. 
The normalization $G(A_0)$ is correct too, since $G_k(A_0)=1$ by construction.
 
We are almost finished now. We started from a sequence $\{ G_k \} = \{ G_k^\infty \}$, 
then extracted a subsequence so that, in particular, $G_k^\infty$ converges to some limit $G$, 
and then proved that $G$ is the desired Green function. 
The same thing would happen if we started from any other subsequence, 
and we would always get the same limit. Since we can always extract convergent 
subsequences, this means that the limit existed already for our initial sequence, 
and is the desired Green function. 

It remains only to show that $G_k$ converges (strongly) in $W^{1,2}_{loc} (\Omega_\infty).$

\begin{lemma} \label{l2a2}
Keep the notations and assumptions of Theorem \ref{l2a4}.
Then for each compact subset $K$ of $\Omega_\infty$,
\begin{equation} \label{2a11}
\lim_{k \to +\infty} \int_K |\nabla G_k - \nabla G_\infty|^2 = 0.
\end{equation}
\end{lemma}

We do not need to worry about putting our weight $w(X)$ in the integral or not, 
because it is bounded from above and below on any compact subset of $\Omega$. 
Similarly, as long as we stay in a compact subset $K$ of $\Omega_\infty$,
our degenerate elliptic operators coincide on $K$ with multiples of a standard elliptic operator, 
with bounds on the ellipticity that depend only on $K$ and coefficients that in any case satisfy \eqref{2a10bis} -- see the discussion near \eqref{eqConv}. That is, Case 1 and Case 2 are identical as far as the estimates strictly inside the domain are concerned, and we will treat them as such.

\noindent {\it Proof of Lemma~\ref{l2a2}.}
Of course it enough to prove \eqref{2a11} when $K$ is the 
closure of a ball $B = B(x,r)$ such that $4B \subset \Omega_\infty$. The functions $G_k$ and $G$ are bounded in $W^{1,2}(3B)$ uniformly in $k$ (for $k$ sufficiently large) by a constant which we will denote by $C(B)$ (see, in particular, \eqref{2c12b} and \eqref{2c40}.

Let us show that $\{G_k \}$ is a Cauchy sequence in $W^{1,2}(B)$.
For this we fix $k$ and $l$ large, and we want to estimate $\nabla (G_k-G_l)$.
We will introduce an intermediate function $U_{kl}$, which coincides with $G_l$
on the sphere $S_\rho=\partial B_\rho = \partial B(x, \rho)$ for some $\rho \in (r, 2r)$
but satisfies the same equation $L_kU_{kl}=0$ as $G_k$ in $B_\rho$.

First we claim that for almost every radius $\rho \in (r,2r)$,
the restriction $g_k$ of $G_k$ to $S_\rho = \d B(x,\rho)$ lies in the space 
$W^{1,2}(S_\rho, d\sigma)$ of functions of $L^2(S_\rho, d\sigma)$, 
which have a distribution gradient $\nabla_T g_k$ 
in $L^2(S_\rho, d\sigma)$, and in 
addition, that the gradient $\nabla_T g_k$ of $g_k$ is given, almost everywhere on 
$S_\rho$, by the restriction of $\nabla G_k$ to $S_\rho$.
Here $\sigma$ denotes the surface measure on $\sigma$, and before we start
we know that the restriction of $\nabla G_k$ is defined $\sigma$-almost-everywhere 
on almost every $S_\rho$ by Fubini. There is no real doubt as to what we mean by restriction,
because $G_k$ is H\"older continuous on $2B$, but in general we could use Fubini to 
say that $g_k$ is defined almost-everywhere on almost every sphere $S_\rho$, by Fubini.

The claim is not hard. We can use spherical coordinates to reduce matters to the case of 
a function in $W^{1,2}(R)$ for some parallelepiped $R$, and use the classical fact that 
then the restriction to almost every hyperplane $P$ parallel to the axes lies in 
$W^{1,2}(P \cap R)$, with partial derivatives given by the restriction to $P$ 
of the partial derivatives of the function. So we skip the details, 
but refer to the proof of Corollary 14.28 in \cite{MSbook} for a similar computation.
We choose $\rho \in  (r,2r)$ so that in addition to the property above
and its analogue for the restriction $g_l$ of $G_l$ to $S_\rho$, we have that
\begin{equation} \label{2c24}
\begin{aligned}
\int_{S_\rho} (|\nabla_T g_k|^2 + |\nabla_T g_l|^2) d\sigma
&\leq \int_{S_\rho} (|\nabla G_k|^2 + |\nabla G_l|^2) d\sigma
\cr&\leq C r^{-1} \int_{2B} |\nabla G_k|^2 + |\nabla G_l|^2 \leq C(B).
\end{aligned}
\end{equation}
See the beginning of the proof of the Lemma where we discuss the uniform bound $C(B)$ 
and recall that while we do not care how $C(B)$ depends on $B$, it does not depend on
$k$ or $l$.

Going further, denote by $W_\rho = \dot W^{1,2}(B_\rho)$ the space of functions 
$F \in L^2(B_\rho)$, with $\nabla F \in L^2(B_\rho)$, equipped with the homogeneous norm  
$||F||_\rho = \big(\int_{B_\rho} |\nabla F|^2\big)^{1/2}$ for $F \in W_\rho$.
We also use the (homogeneous) space $H_\rho = \dot H^{1/2}(S_\rho)$ of functions 
$f\in L^2(S_\rho,d\sigma)$ which have half a derivative in $L^2(S_\rho)$, once again, 
equipped with the corresponding homogeneous norm. 
We need the following classical facts about $H_\rho$. 
First, every $F \in W_\rho$ has a trace $Tr(F)$ in $H_\rho$,
with $||Tr(F)||_{H_\rho} \leq C ||F||_\rho$. In the other direction, every $f \in H_\rho$ has an extension 
$F \in W_\rho$, with $||F||_\rho \leq C ||f||_{H_\rho}$ and $Tr(F) = f$. 
Using these results and the Lax-Milgram's theorem, we know that 
for each elliptic operator $L = -\dv A \nabla$, 
each $f \in H_\rho$ has a unique $L$-harmonic extension $F$ to $B_\rho$, 
i.e., a function $F \in W_\rho$ such that $LF = 0$ (weakly as in \eqref{2a7}) and $Tr(F) = f$. 
In addition, $||F||_\rho \leq C ||f||_{H_\rho}$ (where now $C$ depends also on $r$ and 
the ellipticity constants for $A$).

Interpolating between $L^2(S_\rho)$ and the space of functions with 
tangential derivatives in $L^2(S_\rho)$, we deduce that
\begin{equation} \label{2c25} 
||f||_{H_\rho}^2 
\leq C ||f||_{L^2(S_\rho)} ||\nabla_T f||_{L^2(S_\rho)}.
\end{equation}
This is a rather simple result of  
interpolation between  
Sobolev spaces, but since we are on a ball, we can provide an even more direct argument. Because we are working with a ball, we can also give a description of $H_\rho$ and its norm in terms 
of spherical harmonics. If $f = \sum_j f_j$ is the decomposition of $f$ into spherical harmonics,
where $f_j$ is the part that comes from harmonics of degree $j$,
(and we do not need to further decompose $f_j$ into polynomials),
we know that $||f||_{L^2(S_\rho)}^2 = \sum_j ||f_j||_2^2$, while 
$||\nabla_T f||_{L^2(S_\rho)}^2 = C \sum_j j^2 ||f_j||_2^2$, and the $H_\rho$-norm 
of $f$ is given by $||f\|_{H_\rho}^2 = \sum_j j ||f_j||_2^2$.
Again, we do not need to know how the constants depend on $\rho$, as 
long as they are uniform in $\rho \in (r,2r)$. Thus by Cauchy-Schwarz, \eqref{2c25} holds. 

Let us apply this to the operator $L_k$ and the restriction $g_l$ of $G_l$ to $S_\rho$;
we obtain a function $U_{kl} \in W_\rho$ such that
\begin{equation} \label{2b25}
L_k U_{kl} = 0 \text{ on $B_\rho$, } \  Tr(U_{kl}) = g_l,
\end{equation}
and
\begin{equation} \label{2b26}
||U_{kl}||_\rho \leq C ||g_l||_{H_\rho} \leq C ||\nabla G_l||_{L^2(S_\rho)} \leq C(B)
\end{equation}
because \eqref{2c25} holds, $G_l$  
is continuous and bounded on $2B$, hence $||g_l||_{L^2(S_\rho)} \leq C$,  
and then by definition of $\rho$ and \eqref{2c24}.
We said that $g_l$ is the restriction of $G_l$ to $S_\rho$, and this makes
sense because $G_k$ is continuous, but it is also the trace of $G$, in the sense 
of the operator above (recall that $G_l \in W_\rho$), and the same is true for $g_k$.

Then $U_{kl}-G_k$ is $L_k$-harmonic, lies in the space $W^{1,2}(B)$, 
and the trace of $U_{kl}$ is  $f = g_l -g_k$. Hence,
\begin{equation} \label{2b27}
\begin{aligned}
||U_{kl}-G_k||_\rho &\leq C ||f||_{H_\rho} 
\leq C ||f||_{L^2(S_\rho)}^{1/2} \, ||\nabla_T (g_l-g_k)||_{L^2(S_\rho)}^{1/2}
\cr&\leq C ||G_l-G_k||_{L^\infty(2B)}^{1/2}\, ||\nabla_T (g_l-g_k)||_{L^2(S_\rho)}^{1/2}
\leq C ||G_l-G_k||_{L^\infty(2B)}^{1/2}
\end{aligned}
\end{equation}
as in \eqref{2b26}, \eqref{2c25}. 
We used the fact that $G_l-G_k$ is continuous on $2B$, so 
$||f||_{L^2(S_\rho)} \leq C ||G_l-G_k||_{L^\infty(2B)}$, and then  
\eqref{2c24}.
This will be sufficient to ensure that $||U_{kl}-G_k||_\rho$ tends to zero
since $||G_l-G_k||_{L^\infty(2B)}$ tends to $0$: $G_k$ tends to 
$G_\infty$ uniformly on compact subsets of $\Omega_\infty$.

We now take care of $||G_l-U_{kl}||_\rho$, with an argument that comes from 
the calculus of variations, but which we had to modify  
because maybe the $L_k$ are not symmetric and the solutions of $L_k U_{kl} = 0$ 
do not minimize an energy. 
We first observe that the trace of $G_l-U_{kl}$ on $S_\rho$ vanishes, by definition of $U$ and 
because the trace of $G_l$ is $g_l$. Now we claim that
\begin{equation} \label{2b28}
\int_{B_\rho} \langle A_l \nabla G_l, \nabla(G_l-U_{kl}) \rangle = 0.
\end{equation}
If $G_l-U_{kl}$ were a test function with compact support in $B_\rho$,
this would just be the definition that $G_l$ is a weak solution of $L_l$, as in \eqref{2a7}.
The fact that \eqref{2b28} also holds now follows from the fact that the test functions
with compact support in $B_\rho$ are dense in the subspace $W_{\rho,0} \subset W_\rho$ 
of functions with a vanishing trace. The same reasoning also yields that
\begin{equation} \label{2b29}
\int_{B_\rho} \langle A_k \nabla U_{kl}, \nabla(G_l-U_{kl}) \rangle = 0,
\end{equation}
because $L_k U_{kl} = 0$ weakly in $B_\rho$. Now 
\begin{eqnarray} \label{2b30}
\int_{B_\rho} |\nabla(G_l-U_{kl})|^2 
&\leq& C \int_{B_\rho} \langle A_k \nabla (G_l-U_{kl}), \nabla(G_l-U_{kl}) \rangle
= C \int_{B_\rho} \langle A_k \nabla G_l, \nabla(G_l-U_{kl}) \rangle
\nn\\ & = & C \int_{B_\rho} \langle (A_k - A_l) \nabla G_l, \nabla(G_l-U_{kl}) \rangle
\end{eqnarray}
because $A_k$ is elliptic, and then by \eqref{2b29} and \eqref{2b28}.
We apply Cauchy-Schwarz, simplify, and get that
\begin{equation} \label{2c32}
||G_l-U_{kl}||_\rho \leq C \Big\{\int_{B_\rho} |A_k - A_l|^2 |\nabla G_l|^2 \Big\}^{1/2}.
\end{equation}

If we had  \eqref{2a10}, we could conclude that the expression above tends to zero simply 
pulling out the $L^\infty(B)$ norm of $A_k-A_l$ and using that $G_l$ are uniformly 
bounded in $W^{1,2}(B)$.

Under the weaker assumption \eqref{2a10bis}, we use the reverse $L^p$ inequality for the gradients of solutions. Indeed, there exists $p>2$,  that depends on $n$ and the ellipticity constants for the $L_l$, such that 
\begin{equation} \label{2a12b}
\Big\{\int_{B_\rho} |\nabla G_l|^p \Big\}^{2/p}
\leq \int_{2B} |\nabla G_l|^2 \leq C(B).
\end{equation}
In particular, the  
power $p$ and a constant $C(B)$ do not depend on $l$ 
(see, e.g., \cite{Ken}, Lemma~1.1.12). This, together with  \eqref{2a10bis}, 
shows that \eqref{2c32} tends to zero:
\begin{equation} \label{2c35}
\begin{aligned}
||G_l-U_{kl}||_\rho^2 
\leq C \int_{B_\rho} |A_k - A_l|^2 |\nabla G_l|^2 
&\leq C \Big\{\int_{B_\rho} |A_k - A_l|^{2q} \Big\}^{1/q} 
\Big\{\int_{B_\rho} |\nabla G_l|^p \Big\}^{1/2p}
\cr&\leq C \Big\{\int_{B_\rho} |A_k - A_l|^{2q}  \Big\}^{1/q},
\end{aligned}
\end{equation}
where $q$ is the  
dual exponent of $p/2$. 
The last expression above tends to $0$ because $A_k-A_l$ tends to $0$ in $L^{1}_{loc}$
and is bounded. The rest of the proof stays as before, so
Lemma \ref{l2a2} follows. This finishes the proof of Theorem~\ref{l2a4}.
\qed

\section{Prevalent approximation of $G^\infty$ when $E$ is uniformly rectifiable}
\label{S3}

The assumptions for the next theorem concern the classical situation, i.e., Case 1 above. 
We are given an unbounded domain $\Omega \subset \R^n$, assume that 
$E = \d\Omega$ is Ahlfors regular of dimension $n-1$ (see the definition \eqref{1a1})
and uniformly rectifiable (Definition \ref{d2}), and that $\Omega$ has interior 
corkscrew points and Harnack chains, as in \eqref{2b4}. We assume, furthermore, that the operator is close to a constant coefficient one, in the sense of \eqref{3b3}.

\begin{theorem} \label{t3a1} Let $\Omega\subset \RR^n$ be an unbounded domain 
with an $n-1$ dimensional uniformly rectifiable boundary $E=\po$, and satisfying the interior corkscrew and Harnack chain conditions \eqref{2b4}.
Assume that $L$ is an elliptic divergence form operator 
that satisfies \eqref{3b3}, and denote by $G = G^\infty$ the Green 
function for $L$, with pole at $\infty$. Then $G$ is prevalently close to the distance to a plane 
(as in Definition \ref{d3}) and for every choice of $\beta > 0$ and any AR measure
$\mu$ on $E$, $G$ is prevalently close to $D_{\beta,\mu}$.
\end{theorem}

In this statement $G$ is only determined modulo a multiplicative constant, 
but this does not matter because the conditions \eqref{1a10} and \eqref{1a12} do not see this constant.

\begin{remark} \label{r3b1}
Something like \eqref{3b3} is needed for this statement. Otherwise consider the falsely 
$2$-dimensional example where $\Omega = \{(x,y) \in \R^2 \, ; \,  y > 0 \} \subset \R^2$ 
and $A(x,y) = a(y) I$ for $(x,y) \in \Omega$, where 
$I$ denotes the identity matrix and for instance
$a(y) = 1$ when $y \in [2^{2k},2^{2k+1})$ and $a(y) = 2$ when 
$y \in [2^{2k+1},2^{2k+2})$, $k \in{\mathbb {Z}}$.
The Green function for $L$ is easy to compute: take $G(x,y) = Cg(y)$, where
$g(0)=0$ and $g'(y) = a(y)^{-1}$; this is clearly a solution of $LG=0$, and
the uniqueness of $G^\infty$ does the rest. We can see that $G$ is not prevalently 
close to a distance function because $g$ spends its time oscillating between two affine functions. 
Yet notice that for this example the harmonic measure $\omega^\infty$ is (proportional to) 
the Lebesgue measure on $E$. 

Moreover, there are more complicated examples \cite{MM, P} which show that  
without the Dahlberg-Kenig-Pipher condition above, the absolute continuity of the harmonic 
measure with respect to the Hausdorff measure could fail even in the half-space. 

For all these reasons our conditions on the coefficients are morally optimal. 
\end{remark}

\smallskip\noindent
{\it Proof of Theorem~\ref{t3a1}.} {\bf Part I: comparison to $\dist(X,P)$.}
Let us define some good sets. 
Let $M \geq 1$ and $\varepsilon > 0$ be given, and consider
the set $\cG(\varepsilon, M)$ of all pairs $(x,r) \in E \times (0,+\infty)$ such that all the following 
properties are satisfied. First, there is a hyperplane $P = P(x,r)$ through $x$, such that
\begin{equation} \label{3a5}
d_{x,2Mr}(E, P) \leq \varepsilon 
\end{equation}
(where $d_{x,r}$ is as in \eqref{1a8}), and also, for one of the two half spaces $H$ bounded by $P$,
\begin{equation} \label{3a6}
d_{x,Mr}(\Omega, H) \leq 2\varepsilon. 
\end{equation}
We claim that each $\cG(\varepsilon, M)$ is a Carleson-prevalent set, with a constant 
that may depend (wildly) on $\varepsilon$, $M$, and the various constants implicit in our 
geometric assumptions.

First consider the pairs $(x,r)$ such that we cannot find a plane $P$ as in \eqref{3a5}. 
If $(x,r)$ is such a pair, it is easy to see that $(x,2Mr) \notin \cG_{ur}(2M\varepsilon)$, 
where $\cG_{ur}$ is the good set of Definition \ref{d2}. 
This set of bad pairs satisfies a Carleson packing condition, by Definition~\ref{d2} and 
because it is easy to see that
$\big\{ (x,r) \, ; \, (x, Ar) \in \cB \big\}$ satisfies a Carleson packing condition when $\cB$ does.

Now we claim that if $P$ satisfies \eqref{3a5}, and if $\varepsilon$ is small enough and 
$M$ is large enough (which we may assume) then \eqref{3a6} automatically holds for one 
of the two half spaces bounded by $P$.
That is, the case when $\Omega \cap B(x,Mr)$ is nearly empty, or on the contrary very close 
to the full $B(x,r)$, is excluded by our NTA assumptions. 
In fact, since we assume that $\Omega$ is a 1-sided NTA domain with a uniformly rectifiable boundary, it possesses exterior corkscrew points as well.
Then, on one hand, we can find a corkscrew point 
$Z$ for $B(x,r)$, this point lies far from $E$, hence far from $P$ too, then we select 
the half space $H$ bounded by $P$ that contains $Z$, and already \eqref{3a5} says 
that every point $Y \in H \cap B(x,Mr)$ such that $\dist(Y,P) > \varepsilon Mr$
lies in $\Omega$, because the line segment from $Y$ to $Z$
does not meet $E$. Similarly, all the points 
$Y \in B(x,Mr) \sm H$ such that $\dist(Y,P) > \varepsilon Mr$
lie in the same component 
of $\R^n \sm E$ as $Z'$, the image of $Z$ by the reflection across $P$. 
We just need to see that this component is not $\Omega$, or in other words that 
$Z' \notin \Omega$.
But if $Z' \in \Omega$, there is a Harnack chain from $Z$ to $Z'$, and by definition 
of a Harnack chain it stays inside of $B(x,Mr)$ if $M$ is large enough. 
One of the balls $B_j$ of the chain meets $P \cap B(x,Mr)$, and this is impossible 
if $\varepsilon$ is small enough because $2B_j \subset \Omega$ and 
$\diam B_j \geq C^{-1} \min(\dist(Z,E), \dist(Z',E)) \geq C^{-1} r$. 
So \eqref{3a6} comes for free.

\begin{lemma} \label{l3a2}
For each choice of $\varepsilon > 0$ and $M > 0$, we can find $\varepsilon_1 > 0$, 
$M_1 \geq 1$, $\kappa > 0$, and $K \geq 1$, depending on $n$, $\varepsilon$, $M$, 
the AR constant for $E$, 
and the one sided NTA constants from \eqref{2b4}, 
such that if $(x,r) \in \cG(\varepsilon_1,M_1) \cap \cG_{cc}(\tau,K)$,
then $(x,r)$ lies in the good set $\cG_{Gd}(\varepsilon,M)$ of Definition \ref{d3}.
\end{lemma}

\noindent {\it{Proof of Lemma~\ref{l3a2}}}. 
Let the various constants and $\varepsilon$, $M$, be given, and suppose the lemma fails 
for these constants. Then there exists an example $(\Omega_k, E_k, L_k)$ where all 
the assumptions are satisfied, and for some pair $(x_k,r_k) \in E_k \times (0,+\infty)$, 
we have that
$(x_k,r_k) \in \cG(\varepsilon_k,M_k) \cap \cG_{cc}(\tau_k,K_k)$, say, with
$\varepsilon_k = \tau_k = 2^{-k}$ and $K_k = M_k = 2^k$, and yet
$(x_k,r_k) \notin \cG_{Gd}(\varepsilon,M)$. 

As the reader guessed, we want to take a limit and derive a contradiction.
By translation, dilation, and rotation invariance of the problem, we may assume that
$x_k = 0$, $r_k = 1$, and we can take a fixed hyperplane $P$ through the origin and a fixed
half space $H$ bounded by $P$ that work in \eqref{3a5} and \eqref{3a6} for all $k$.
By \eqref{3a5} and \eqref{3a6}, $E_k$ tends to $P$ and $\Omega_k$ tends to $H$, as in \eqref{2a2}.
Pick a corkscrew point $X_0$, for instance the point of $H$ that lies at distance $1$ from both $0$ and $P$.
Then let $G_k$ denote the Green function for $\Omega_k$, with the pole at $\infty$, and normalized
by $G_k(X_0)=1$.
We want to apply Theorem \ref{l2a4}, so we check the last hypothesis \eqref{2a10bis}.
The assumption that $(x_k,r_k) \in \cG_{cc}(\tau_k,K_k)$ gives us a constant matrix $A_{k,0}$ 
and, modulo extracting a subsequence, we may assume that $A_{k,0}$ tends to a limit $A_{0,0}$.
Now \eqref{2a10bis} holds, with $A_\infty = A_{0,0}$, because if $B$ is a ball such that 
$2B \subset \Omega_\infty$, then for $k$ large $B \subset W_{K_k}(0,1) = W_{K_k}(x_k,r_k)$,
and 
$$\int_{B} |A_k-A_\infty| \leq |B| |A_{k,0}-A_{0,0}| + \int_B |A_{k}-A_{k,0}|
\leq |B| |A_{0,k}-A_{0,0}| + \tau_k, 
$$
which tends to $0$.

So Theorem \ref{l2a4} says that 
$\{ G_k \}$ converges, uniformly on compact subsets of $\RR^n$,
to the Green function $G_\infty$ for $L_\infty = - \dv A_{0,0} \nabla$
on $H = \Omega_\infty$, with a pole at $\infty$ and normalized at $X_0$. 
Now $G(X) = \dist(X,P)$ is a solution to $L_\infty G = 0$, 
lies in all the correct $W_r(B)$, vanishes along $P$, and has the same normalization, so 
$G_\infty(X) = \dist(X,P)$ on $\Omega_\infty = H$. 
Thus the $G_k$ converge uniformly, on the compact subsets of $\RR^n$,
to the distance function. Therefore, for $k$ large,
\begin{equation} \label{3a10}
|\dist(X,P) - G_k(X)| \leq \varepsilon
\text{ for } X \in \Omega \cap B(0,M),
\end{equation}
which contradicts the fact that \eqref{1a10} fails for $G_k$ and the unit ball.
This completes the proof of Lemma \ref{l3a2}. \qed

\smallskip
The first part of Theorem \ref{t3a1}, with the distance function to a plane, follows at once, 
by a combinaison of definitions and because the union 
$\cG(\varepsilon_1,M_1) \cap \cG_{cc}(\tau,K)$ of two prevalent sets is prevalent.

\noindent {\it Proof of Theorem~\ref{t3a1}.} {\bf Part II: comparison to $D_{\beta,\mu}(X)$.}
We now intend to show that $G$ is prevalently close to any given $D_{\beta,\mu}$.
To this end, we will 
use the result above, and compare $D_{\beta,\mu}$ to $\dist(X,P)$. 
That is, we will show that prevalently, the restriction of $D_{\beta,\mu}$ to $\Omega$
is close to some $\dist(X,P)$, and then we will deal with organization issues, for instance 
making sure that we use the same plane in the two descriptions.

The estimates that follow are similar to the proofs that were done in \cite{DEM}, 
but they are a little simpler here because we can content ourselves with weak estimates.
We intend to use a result of X. Tolsa on the good approximation of $\mu$ by flat measures,
and for this we need some notation. 
These results will be used later in full generality, and for the time being we
assume that $E$ is an Ahlfors regular set $E$ dimension $d$, with any
$d \in (0,n)$, and $\mu$ is any AR measure of dimension $d$ whose support is $E$.

Denote by $\cF_d$ the set of flat measures of dimension $d$ in $\R^n$, 
i.e., measures $\sigma = c \H^d_{\vert P}$, where $c>0$ and $P$ is an affine $d$-plane.
For $x \in \R^n$ and $r > 0$, denote by $\Lambda(x,r)$ the set of Lipschitz functions
$\varphi : \R^n \to \R$ such that $\varphi(x)=0$ on $\R^n \sm B(x,r)$ 
and the Lipschitz norm of $\vp$ is less than or equal to $1$. 
Then define a sort of Wasserstein distance between two measures $\mu$ and $\sigma$ by
\begin{equation} \label{3a12}
D_{x,r}(\mu, \sigma) = r^{-d-1} \sup_{\varphi \in \Lambda(x,r)} 
\Big| \int \varphi d\mu - \int \varphi d\sigma\Big|,
\end{equation} 
and finally set
\begin{equation} \label{3a13}
\alpha(x,r) = \inf_{\sigma \in \cF_d} D_{x,r}(\mu, \sigma).
\end{equation}
The normalization is such that when $\mu$ is an AR measure, 
$\alpha(x,r) \leq C$ systematically, but when $E$ is uniformly rectifiable and 
$\mu$ is any AR measure on $E$, $\alpha(x,r)$ is often much smaller than this.
In fact, it is proved in \cite{To} 
that $\alpha^2(x,r) \frac{d\mu(x) dr}{r}$ is a Carleson measure on $E \times (0,+\infty)$.
Here we shall only use the corresponding weak estimate, which follows by Chebyshev's inequality, 
that says that for each choice of $\eta > 0$ and $N \geq 1$, 
\begin{equation} \label{3a14}
\cG_{\alpha}(\eta,N) = \big\{ (x,r) \in E \times (0,+\infty) \, ; \, \alpha(x,Nr) \leq \eta \big\}
\ \text{ is a Carleson-prevalent set.}
\end{equation}

\begin{lemma} 
\label{l3a14}
Let $E$ and $\mu$ be Ahlfors regular of dimension $d\in (0, n)$ in $\R^n$. Let
$(x,r) \in \cG_{\alpha}(\eta,N) \subset E \times (0,+\infty)$, and let 
$\sigma = c \H^d_{\vert P}$ be a flat measure such that
\begin{equation} \label{3d15}
D_{x,Nr}(\mu, \sigma) \leq 2\eta.
\end{equation}
If $N$ is large enough (depending on $n$, $d$, and the AR 
constant for $\mu$) and $\eta$ is small enough (depending on the aforementioned 
parameters and also on $N$), then
\begin{equation} \label{3a17}
\dist(z,P) \leq \eta_1 r \ \text{ for $z \in E \cap B(x,Nr/3)$,} 
\end{equation}
and 
\begin{equation} \label{3a17bis}
\dist(z,E) \leq \eta_1 r \ \text{ for $z \in P \cap B(x,Nr/3)$,} 
\end{equation}
where we set $\eta_1 = C_1 N \eta^{1/(d+1)}$, with a constant $C_1$ that depends 
on the AR constant for $\mu$, $d$ and $n$.  Furthermore,
for the distance function $D_{\beta,\mu}$,
\begin{equation} \label{3d17}
|D_{\beta,\mu}(X) - (ca_\beta)^{-1/\beta}\dist(X,P)| \leq C_2 N \eta^{1/(d+2+\beta)} r
+C_2 \dist(X,P)^{1+\beta} N^{-\beta} r^{-\beta}
\end{equation}
for $z \in B(x,Nr/5)$, where $C_2$ depends on $n$, $d$, the AR constant for $\mu$, 
and $\beta$, and $a_\beta$ is a geometric constant that depends on $\beta$.
\end{lemma}

\bp
The definition of $D_{x,Nr}(\mu, \sigma)$ says that
\begin{equation} \label{3d18}
\Big| \int \varphi \big(d\mu - d\sigma \big) \Big| 
\leq (Nr)^{d+1}D_{x,r}(\mu, \sigma) \leq 2 \eta N^{d+1} r^{d+1}
\end{equation}
for every $\varphi \in \Lambda(x,Nr)$. We intend to use this for various functions $\varphi$
to get relevant information on $\mu$ and $\sigma$. Let us first check that 
\begin{equation} \label{3d20}
\text{$P$ meets $B(x,Nr/4)$ and $C^{-1} \leq c \leq C$.}
\end{equation}
Here and throughout the proof,  
$C$ is a constant which depends on $d,n$, and the AR constants of $\mu$ only, 
and whose value may change  
from line to line.
Take a first bump function $\psi_1$ such that $0 \leq \psi_1 \leq 1$ everywhere,
$\psi_1(z) = 0$ on the complement of $B(x,Nr/4)$,  
and $\psi_1(z) = 1$ on $B(x,Nr/8)$.
We can manage to do this with $|\nabla \psi| \leq (Nr/8)^{-1}$, so
$(Nr/8) \psi_1 \Lambda(x,Nr)$ and \eqref{3d18} says that
\begin{equation} \label{3d19}
\Big| \int \psi_1 \big(d\mu - d\sigma \big) \Big| \leq 16 \eta N^{d} r^{d}.
\end{equation}
Since $\psi_1(z) = 1$ on $B(x,Nr/8)$ and $\mu(B(x,Nr/8) \geq C^{-1} N^d r^d$,
we see that (if $\eta$ is small enough) $\sigma(B(x,Nr/4)) \geq (2C)^{-1} N^d r^d$
(because $\psi_1(z) = 0$ on $\R^n \sm B(x,Nr/4)$), and hence
$P$ meets $B(x,Nr/4)$ and $c \geq C^{-1}$.

Similarly, if now $\psi$ is a similar function, but with $0 \leq \psi \leq 1$ everywhere,
$\psi(z) = 0$ on $\R^n \sm B(x,Nr)$, $\psi(z) = 1$ on $B(x,Nr/2)$, and 
$|\nabla \psi| \leq (Nr/2)^{-1}$, we get that
\begin{equation} \label{3d21}
\Big| \int_{B(x,Nr)} \psi \big(d\mu - d\sigma \big) \Big| \leq 4 \eta N^{d} r^{d}.
\end{equation}
But $\mu(B(x,Nr)) \leq C N^d r^d$, so $\int_{B(x,Nr)} \psi d\sigma \leq 2C N^d r^d$,
and since $\int_{B(x,Nr)} \psi d\sigma \geq \sigma(B(x,Nr/2)) = c \H^d(P \cap B(x,Nr/2))$
and $P$ meets $B(x,Nr/4)$, we get that $c \leq C$ too; \eqref{3d20} follows.

Next we want to control the distance to $P$. 
We keep the same function $\psi$, but try the product $\varphi(z) = \psi(z) \dist(z,P)$. 
We still have that $\varphi = 0$ on $\R^n \sm B(x,Nr)$, as required in the definition 
of $\Lambda(x,Nr)$, but now
\begin{equation} \label{3c}
|\nabla \varphi(z)| \leq  \psi(z) +  |\nabla \psi(z)| \dist(z,P) 
\leq 1 + (Nr/2)^{-1} \dist(z,P) \leq 5
\end{equation}
because $P$ meets $B(x,Nr/4)$. So $\varphi/5 \in \Lambda(x,Nr)$, and \eqref{3d18} yields
\begin{equation} \label{3c23}
\Big| \int \psi(z) \dist(z,P) \big(d\mu(z) - d\sigma(z) \big) \Big| \leq 10 \eta N^{d+1} r^{d+1}.
\end{equation}
Now $\int \psi(z) \dist(z,P)  d\sigma = 0$ because $\sigma$ is supported on $P$,
and since $\psi(z) = 1$ on $B(x,Nr/2)$, we get that
\begin{equation} \label{3c24}
\int_{B(x,Nr/2)}\dist(z,P) d\mu(z) \leq 10 \eta N^{d+1} r^{d+1}.
\end{equation}
We still need to transform this into an $L^\infty$ bound. So we assume
that \eqref{3a17} fails, so that $\dist(z,P) \geq \eta_1 r$ for some
$z\in E \cap B(x,Nr/3)$. Then $\dist(z,P) \geq \eta_1 r/2$ on $B_1 = B(z,\eta_1 r/2)$,
and since $\eta_1 < Nr/6$ if $\eta$ is small enough, $B_1 \subset B(x,Nr/2)$ and
\eqref{3c24} yields
\begin{equation} \label{3c25}
10 \eta N^{d+1} r^{d+1} \geq (\eta_1 r/2) \mu(B_1) \geq C^{-1} (\eta_1 r)^{d+1}
= C^{-1} (C_1 N \eta^{1/(d+1)})^{d+1} r^{d+1}
\end{equation}
by Ahlfors regularity and the definition of $\eta_1$.
We now choose $C_1$ so large that \eqref{3c25} is impossible,
and this contradiction completes the proof of \eqref{3a17}.

For \eqref{3a17bis} we proceed exactly the same way, but with the function 
$\psi(z) \dist(z,E)$, which of course vanishes on the support of $\mu$.
Maybe we need to make $C_1$ a little larger because we now use the Ahlfors regularity
of $\sigma$, but this is all right.

\ms
Now we estimate $D_{\beta,\mu}(X) - (ca_\beta)^{-1/\beta} \dist(X,P)$. 
Recall that $D_{\beta,\mu} = D_\beta(X) = R_\beta(X)^{-1/\beta}$, with 
\begin{equation} \label{3c26}
R_\beta(X) = \int_E |X-y|^{-d-\beta} d\mu(y),
\end{equation}
so we take care of $R_\beta$ first. Let $R(X) = \int_E |X-y|^{-d-\beta} d\sigma(y)$
denote the analogue of $R_\beta(X)$, but for the measure $\sigma$. 
A direct computation shows that
\begin{equation} \label{3c27}
R(X) = c a_\beta \dist(X,P)^{-\beta},
\end{equation}
just because $\sigma = c \H^d_{\vert P}$. Let  
$\psi$ be a smooth cut-off function such that $0 \leq \psi \leq 1$ everywhere,
$\psi = 0$ outside of $B(x,Nr/4)$, $\psi = 1$ on $B(x,Nr/8)$,
$|\nabla \psi| \leq 16(Nr)^{-1}$ and $|\nabla^2 \psi| \leq C(Nr)^{-2}$ 
(to be used later). We first study the main part of $R_\beta(X) - R(X)$, namely
\begin{equation} \label{3c28}
A(X) = \int |X-y|^{-d-\beta} \psi(y) [d\mu(y)-d\sigma(y)].
\end{equation}
Also let $\eta_2 > \eta_1$ be small, to be chosen later, and then introduce a 
smooth cut-off function $\xi$ such that $0 \leq \xi \leq 1$ on $\R^n$, 
$\xi(X) = 1$ when $\dist(X,P) \leq \eta_2 r$, $\xi(X) = 0$ when $\dist(X,P) \geq 2\eta_2 r$,
$|\nabla \xi| \leq 2(\eta_2 r)^{-1}$, and $|\nabla^2 \xi| \leq C(\eta_2 r)^{-2}$. 
Notice that $\xi(y) = 1$ when $y \in B(x,Nr/4)$ lies in the support of 
$\sigma$ or $\mu$, by \eqref{3a17}, so 
\begin{equation} \label{3c29}
A(X) = \int |X-y|^{-d-\beta} \psi(y) \xi(y) [d\mu(y)-d\sigma(y)].
\end{equation}
We are going to restrict our attention to the case when
\begin{equation} \label{3c30}
X \in H = \big\{ X \in B(x,Nr/5)\, ; \, \dist(X,P) \geq 3\eta_2 r \big\},
\end{equation}
because in this region the distance from $X$ to the support of $\xi$ is at least
$\eta_2 r$. Set 
$\varphi_X(y) = |X-y|^{-d-\beta} \psi(y) \xi(y)$ for $X \in H$. Notice that
$\varphi_X$ is supported on $B(x,Nr)$ because of $\psi$, and 
\begin{equation} \label{3c31}
\begin{aligned}
|\nabla \varphi_X(y)| &\leq C |X-y|^{-d-\beta-1} + C (\eta_2 r)^{-1} |X-y|^{-d-\beta}  
+16 (Nr)^{-1} |X-y|^{-d-\beta}
\cr&\leq C (\eta_2 r)^{-1} |X-y|^{-d-\beta} \leq C (\eta_2 r)^{-d-\beta-1}
\end{aligned}
\end{equation}
because $|X-y| \geq \eta_2 r$. Normally we should cut the integral into annuli and get a 
better result, but let us not bother. We apply \eqref{3d18} to a multiple of $\varphi_X$ 
and get that
\begin{equation} \label{3c32}
\begin{aligned}
|A(X)| &= \big|\int \varphi_X [d\mu-d\sigma]\big| 
\leq 2 \eta (Nr)^{d+1} ||\nabla \varphi_X||_\infty
\cr&\leq C \eta (Nr)^{d+1} (\eta_2 r)^{-d-\beta-1}
= C \eta N^{d+1} \eta_2^{-d-\beta-1} r^{-\beta}.
\end{aligned}
\end{equation}
This will   
be good enough (because we can choose $\eta$ last). For the rest of the integral,
set $B(X) = R_\beta(X) - R(X) - A(X)$, and notice that
\begin{equation} \label{3c33}
\begin{aligned}
|B(X)| &= \Big|\int |X-y|^{-d-\beta} (1-\psi(y)) [d\mu(y)-d\sigma(y)] \Big|
\cr&\leq \int_{\R^n \sm B(x,Nr/8)} |X-y|^{-d-\beta} (1-\psi(y)) [d\mu(y)+d\sigma(y)]
\leq C (Nr)^{-\beta}
\end{aligned}
\end{equation}
by the Ahlfors regularity of both measures, and a simple argument where one cuts
the domain into annuli of size $2^l Nr$. We add this to \eqref{3c32} and get that
\begin{equation} \label{3c34}
|R_\beta(X) - R(X)| \leq C [\eta N^{d+1} \eta_2^{-d-\beta-1} + N^{-\beta}] r^{-\beta}. 
\end{equation}
Set $D(X) = R(X)^{-1/\beta} = (ca_\beta)^{-1/\beta} \dist(X,P)$ by \eqref{3c27},
and notice that 
\begin{equation} \label{3d35}
D_{\beta,\mu}(X) - D(X) = R_\beta(X)^{-1/\beta} - R(X)^{-1/\beta}.
\end{equation}
Set $t = R_\beta(X)$ and $u=R(X)$ for a minute, and
observe that since the derivative of 
$t \mapsto t^{-1/\beta}$ is $-\frac{1}{\beta} t^{-\frac{1+\beta}{\beta}}$, we 
get that $|t^{-1/\beta}-u^{-1/\beta}| 
\leq C |t-u| (t^{-\frac{1+\beta}{\beta}} + u^{-\frac{1+\beta}{\beta}})$
(including in the case when $|u-t|$ is not small compared to $u$ and $t$).
Here $C^{-1} \dist(X,P)^{-\beta} \leq R(X) \leq C \dist(X,P)^{-\beta}$ 
by direct calculation, and the same thing holds for $R_\beta(X)$, because $\mu$ is 
Ahlfors regular and $\dist(X,P) \sim \dist(X,E)$ for $X\in H$, so
$t$ and $u$ are both comparable to $\dist(X,P)^{-\beta}$, and \eqref{3c34} yields
\begin{equation} \label{3c37}
\begin{aligned}
|D_{\beta,\mu}(X) - D(X)| \leq C \dist(X,P)^{1+\beta}
[\eta N^{d+1} \eta_2^{-d-\beta-1} + N^{-\beta}] r^{-\beta}.
\end{aligned}
\end{equation}
Let us choose $\eta_2 = C_1 N \eta^{1/(d+2+\beta)}$. This way our constraint that 
$\eta_2 \geq \eta_1$ is clearly satisfied (we just decreased the power a little bit),
and \eqref{3c37} simplifies to
\begin{equation} \label{3c38}
\begin{aligned}
|D_{\beta,\mu}(X) - D(X)| 
&\leq C \dist(X,P)^{1+\beta}
[\eta^{1 - \frac{d+1+\beta}{d+2+\beta}} N^{-\beta} + N^{-\beta}] r^{-\beta} 
\cr&\leq C \dist(X,P)^{1+\beta} N^{-\beta} r^{-\beta}.
\end{aligned}
\end{equation}
This is good enough for \eqref{3d17}, but we only get it 
for $X \in H$. When $x \in B(r,Nr/5) \sm H$, we know that
$\dist(X,P) \leq 3 \eta_2 r$ by definition of $H$, and then $\dist(X,E) \leq 4 \eta_2 r$
by \eqref{3a17bis} and because $\eta_2 \geq \eta_1$, hence 
$D_{\beta,\mu}(X) \leq C \dist(X,E) \leq C \eta_2 r$ (because $D_\beta(X)$ is 
always comparable to $\dist(X,E)$). So 
\begin{equation} \label{3c39}
|D_{\beta,\mu}(X) - D(X)| \leq C \eta_2 r = C N \eta^{1/(d+2+\beta)} r
\end{equation}
for $X \in B(r,Nr/5) \sm H$.  Now \eqref{3d17} follows from this and \eqref{3c38},
and this completes our proof of Lemma \ref{l3a14}.
\qed

\begin{remark} \label{r31}
The proof of \eqref{3c34} also shows that 
\begin{equation} \label{3c41}
|\nabla R_\beta(X) - \nabla R(X)| 
\leq C [\eta N^{d+1} \eta_2^{-d-\beta-2} + N^{-\beta-1}] r^{-\beta-1}.
\end{equation}
This is because when we compute $\nabla R$ or $\nabla R_\beta$,
we replace the kernel $|X-y|^{-d-\beta}$ with a kernel of size
$|X-y|^{-d-\beta -1}$; we can then follow the computations above with this 
different power. That is, the estimates are the same as they would be with
$D_{\beta+1}$. Then the proof of \eqref{3c38} yields 
\begin{equation} \label{3c42}
\begin{aligned}
|\nabla D_\beta(X) - \nabla D(X)|  
\leq C \dist(X,P)^{2+\beta} N^{-\beta-1} r^{-\beta-2}.
\end{aligned}
\end{equation}
for $X\in B(x,Nr/5)$ such that $\dist(X,P) \geq C \eta'_2 r$, with now
$\eta'_2 = C_1 N \eta^{1/(d+3+\beta)}$. 
\end{remark}

\ms
We may now return to the proof of prevalence for the good approximation of $G$
by $D_{\beta,\mu}$, that we need to complete our proof of Theorem \ref{t3a1}.
Let $\varepsilon > 0$, $M \geq 1$ and $(x,r) \in E \times (0,+\infty)$ be given.
We want to show that prevalently, $(x,r)$ lies in the good set 
$\cG_{GD_\beta}(\varepsilon, M)$ of Definition~\ref{d3}. 
We may assume that $(x,r) \in \cG_{Gd}(\varepsilon,M)$
because we already proved that this happens prevalently, 
or rather that \eqref{1a10} holds with $\varepsilon/3$. 
We may also assume that $(x,r) \in \cG_{\alpha}(\eta,N)$,
with values of $\eta$ and $N$ that will be chosen soon, because \eqref{3a14} says that
this happens prevalently. Then Lemma \ref{l3a14} says that \eqref{3a17}-\eqref{3d17}
hold, but maybe this happens with a different plane $P_1$. So we will need to see
whether $P_1$ is close enough to $P$.

First observe that in the proof of prevalence for $\cG_{Gd}(\varepsilon/3,M)$,
we can actually find $P$ such that in addition to \eqref{1a10}, we also have that
\begin{equation} \label{3c43}
d_{x, Mr}(E,P) \leq (10M)^{-1} \varepsilon ;  
\end{equation}
see the proof of \eqref{3a5} and \eqref{3a6}, using the definition of uniform rectifiability.
Now \eqref{3a17} and \eqref{3a17bis} also imply that 
\begin{equation} \label{3c44}
d_{x, Mr}(E,P_1) \leq  (N/3) M^{-1} d_{x, Nr/3}(E,P_1) \leq M^{-1} \eta_1
\leq C N M^{-1} \eta^{1/(d+1)} \leq (10M)^{-2} \varepsilon
\end{equation}
because we can pick $N > 3M$, and choose $\eta$ small enough, depending on 
$M$, $N$, and $\varepsilon$. So $P_1$ and $P$ are close to each other, and
\begin{equation} \label{3c45}
|\dist(X,P_1)| - \dist(X,P)| \leq \varepsilon r/3
\end{equation}
for $X\in B(x,M)$. Finally, \eqref{3d17} says that for $X \in B(x,Mr) \subset B(x,Nr/3)$,
\begin{equation} \label{3d46}
\begin{aligned}
|(ca_\beta)^{-1/\beta} D_{\beta,\mu}(X) - \dist(X,P_1)| 
&\leq C N \eta^{1/(d+2+\beta)} r + C \dist(X,P_1)^{1+\beta} N^{-\beta} r^{-\beta}
\cr&\leq C N \eta^{1/(d+2+\beta)} r + C (Mr)^{1+\beta} N^{-\beta} r^{-\beta}
\end{aligned}
\end{equation}
(recall that $c \geq C^{-1}$ by \eqref{3d20}). We choose so $N$ large, 
depending on $M$ and $\beta$, that $CM^{1+\beta} N^{-\beta} < \varepsilon/6$,
and then $\eta$ so small, depending on $M$ and $N$, that 
$C N \eta^{1/(d+2+\beta)} < \varepsilon/6$. Then
$|(ca_\beta)^{-1/\beta} D_{\beta,\mu}(X) - \dist(X,P_1)| \leq  \varepsilon r/3$, and
\eqref{1a12} follows from \eqref{3c45} and \eqref{1a10}.
That is, $(x,r) \in \cG_{GD_\beta}(\varepsilon,M)$.
This completes the proof of Theorem \ref{t3a1}.
\qed

\ms
The next statement is the version of Theorem \ref{t3a1} in Case 2.

\begin{theorem} \label{t3a1bis}
Let $\Omega$ be a domain in $\R^n$, whose boundary $E = \d \Omega$ is
an Ahlfors regular and uniformly rectifiable set of (integer) dimension $d \in (0,n)$. 
If $d = n-1$, assume, in addition, that $E$ is satisfies the interior corkscrew point and 
Harnack chain conditions, as in \eqref{2b4}. 
Let $\alpha > 0$ be given, pick any AR 
measure $\mu$ on $E$ (as in \eqref{1a1}), and define $D_\alpha$ and 
$L = - \dv A D_\alpha^{d+1-n} \nabla$ by
\eqref{1a4}, \eqref{1a3}, and \eqref{2b6-bis}, with $A$ satisfying \eqref{2b2}, \eqref{2b3}, and \eqref{3b3} with $A_0\equiv I$. 
Denote by $G = G^\infty$ the Green function for $L$, with pole at $\infty$. 
Then $G$ is prevalently close to the distance to a plane 
(as in Definition \ref{d3}) and for every $\beta > 0$
and every AR measure $\nu$ on $E$,  
$G$ is prevalently close to  $D_{\beta,\nu}$.
\end{theorem}

As was observed in the comments below Theorem \ref{ti1}, 
the special operators $L$ of \eqref{1a2} are built with specific distance functions $D_\alpha$
that respect the rotation invariance of $\R^n$; this is the reason why we restrict to $A_0 = I$ here, while we authorized other constant coefficient operators in Theorem \ref{t3a1}.

See Section \ref{S2} for the definition and some information on $G^\infty$.
Notice that the measure $\nu$ used to compute the distance $D_{\beta,\nu}$ 
does not need to be the same as the measure $\mu$ used to define $D_\alpha$ and $L$
because the two measures are used in different places and never interact.

\ms

\bp The proof will be nearly the same as for Theorem \ref{t3a1}.
We made sure, when we treated Part II (about $D_\beta = D_{\beta,\nu}$) 
not to use the specific dimension of $E$, and even less the precise form 
of the operator. In particular Lemma \ref{l3a14} is still valid in the present situation, and
this means that when $(x,r) \in \cG_{\alpha}(\eta,N)$ (for the measure $\nu$), 
and $\eta$ and $N$ are chosen correctly, we can find a $d$-plane $P$ such that 
\begin{equation} \label{3c48}
d_{x, M_1 r}(E,P) \leq (10M)^{-2} \varepsilon_1, 
\end{equation}
as in \eqref{3c43}, and where $\varepsilon_1$ and $M_1$ are given in advance.
We also require that $(x,r) \in \cG_{\alpha}(\eta,N)$ for the measure $\mu$
and we claim that (if $\eta$ and $N$ are chosen correctly again), then automatically 
\begin{equation} \label{3c49}
(x,r) \in \cG_{Gd}(\varepsilon,M),
\end{equation}
i.e., we have 
a good approximation of $G$ by $\dist(X,P)$ as required in \eqref{1a10},
in fact with the the same $P$ as in \eqref{3c48}.

We prove this as we did for Lemma \ref{l3a2}, by contradiction and compactness. 
If we cannot get \eqref{3c49}, there exist examples $(\Omega_k, E_k, \nu_k, L_k)$ 
where all the general assumptions are satisfied, 
yet for some pair $(x_k,r_k) \in E_k \times (0,+\infty)$, we have that
\eqref{3c48} holds with $M = 2^k$ and $\varepsilon = 2^{-k}$,
$(x_k,r_k) \in \cG_{\alpha}(2^{-(d+2)k},2^k)$, but
$(x_k,r_k) \notin \cG_{Gd}(\varepsilon,M)$.

We may assume that $x_k = 0$, $r_k = 1$, and that all the planes
$P_k$ coming from \eqref{3c48} are parallel to a same plane $P$
through the origin (we did not require $P_k$ to contain $x_k$).
Then by \eqref{3c48}, $\{ E_k \}$ converges to $P$. If $d=n-1$, we can proceed as 
in the proof of Theorem \ref{t3a1} to extract a subsequence (or turn the domains)
so that $\Omega_k$ converges to a fixed half space $H$ bounded by $P$.
If $d\leq n-2$, a direct inspection shows that $\Omega_k$ converges to
$\R^n \sm P$.

We required that $(x_k,r_k) \in \cG_{\alpha}(2^{-(d+2)k},2^k)$ (for $\nu_k$) 
because this way, we get the existence of a flat measure $\sigma_k$ such that by \eqref{3a12}
\begin{equation} \label{3b50}
\Big| \int \varphi d\nu_k - \int \varphi d\sigma_k\Big|
\leq 2^{(d+1)k} D_{0,2^k}(\nu_k, \sigma_k) \leq 2^{-k}
\end{equation} 
for every $1$-Lipschitz function $\varphi$ supported on $B(0,2^k)$.
By construction (and in particular the proof of \eqref{3c48}), we could
make sure that $\sigma_k$ is supported by $P_k$, but even if we were not that cautious,
it follows from the proof of \eqref{3a17} and \eqref{3a17bis} that the support $P'_k$
of $\sigma_k$ tends to $P$ too. In other words, the limiting measure of $\nu_k$ guaranteed by Theorem~\ref{l2a4} is, in fact, flat.

Similarly to the argument of Theorem~\ref{l2a4} in the paragraph right after \eqref{2c41}, 
it follows easily from \eqref{3b50} (and the uniform convergence of the integrals 
in \eqref{1a3} at $\infty$) that the functions $R_{\alpha, \nu_k}$ of \eqref{1a3} 
(but associated to $\nu_k$) converge, uniformly on every compact subset of $\R^n \sm P$, 
to the function $R_{\alpha, \sigma}$ associated to the limit of the $\sigma_k$ 
(maybe after extraction of a subsequence so that the coefficients $c_k$ converge).

Now look at the operators. We also required that $(x,r) \in \cG_{\alpha}(\eta,N)$ for the measure $\mu$, and so we may assume, by extracting a subsequence again, that
$\mu_k$ converges to a limit $\mu_\infty$ which
is a flat measure on the same plane $P$ as before.  
Then by Theorem~\ref{l2a4} the (properly normalized) Green functions with the pole 
at infinity of $L_k=-\dv D_{\alpha,\mu_k}^{d+1-n} \nabla$ and 
$L_\infty =- \dv D_{\alpha,\mu_\infty}^{d+1-n} \nabla$ converge uniformly
on compact subsets of $\RR^n$, to the Green function with the pole at infinity for 
$L_\infty = -\dv A_k D_{\alpha,\mu_\infty}^{d+1-n} \nabla$ on the domain $\Omega_\infty$. Here $\mu_\infty$ is a flat measure of a plane $P$, $\Omega_\infty$ is $H$ if $d=n-1$ and 
$\R^n \sm P$ otherwise, and so the Green function with the pole at infinity for the limiting operator is $c \dist(X,P)$. This contradicts our assumption that
\begin{equation} \label{3b52}
|c_k G_k(X_k) - \dist(X,P)| \geq \varepsilon
\end{equation}
for some $X_k \in \Omega_k \cap B(0,M)$ and for any $c_k$.  
It is here that we use $A_0=I$ in the approximation of $A_k$ in order to ensure 
that the solution of the limiting equation in the exterior of $\RR^d$ is the distance to the boundary.

So the Green function $G$ is prevalently close to the distance to a plane, 
and now the same argument as for 
Theorem \ref{t3a1} shows that it
is also prevalently close to any $D_{\beta,\nu}$.
Theorem~\ref{t3a1bis} follows.
\qed

\section{Approximation of $\nabla G^\infty$ when $E$ is uniformly rectifiable}
\label{S4}

We want to extend the positive results of Section \ref{S3} to a more precise notion of approximation of $G$ with distance functions, where we also control the first derivative
of $G$. This time, we will not control $\nabla G$ all the way to the boundary, so our
notion of approximation will use the Whitney regions associated to balls $B(x,r)$
centered on the boundary $E$, and defined as in \eqref{3b2} by 
\begin{equation} \label{4a1}
W_M(x,r) = \big\{ X \in \Omega \cap B(x,Mr) \, ; \, \dist(X,E) \geq M^{-1} r \big\}.
\end{equation}
Also, we only control (a little bit more than) the $L^2_{loc}$ norm of $\nabla G$,
so we will use that norm. Except for this, the following definition is similar to Definition \ref{d3} above.

\begin{definition} \label{d41}
Let $G^\infty$ denote the Green function for the operator $L$ in the domain
$\Omega$ bounded by an Ahlfors regular set $E$ of dimension $d < n$.
We say that $\nabla G^\infty$ is prevalently close to the gradient of the distance to a plane 
when for each choice of $\varepsilon >0$ and $M \geq 1$, the set 
$\cG_{\nabla Gd}(\varepsilon, M)$ 
of pairs $(x,r) \in E \times (0,+\infty)$ such that there exists a $d$-plane $P(x,r)$ and a positive constant  $c > 0$, with
\begin{equation} \label{4a3}
\int_{W_M(x,r)} |\nabla \dist(X,P) - c \nabla G^\infty(X)|^2 \leq \varepsilon r^{n},
\end{equation}
is a Carleson-prevalent set.

If in addition we are given an AR measure $\nu$ on $E$ and an exponent $\beta>0$,
and $D_{\beta,\nu}$ is defined as in \eqref{1a3}-\eqref{1a4}, 
we say that $\nabla G^\infty$ is prevalently close to $\nabla D_{\beta,\nu}$ 
when for each choice of $\varepsilon >0$ and $M \geq 1$, the set 
$\cG_{\nabla GD_\beta}(\varepsilon, M)$ 
of pairs $(x,r) \in E \times (0,+\infty)$ such that there exists a positive constant $c > 0$, with
\begin{equation} \label{4a4}
\int_{W_M(x,r)} |\nabla D_{\beta,\nu}(X) - c \nabla G^\infty(X)|^2 dX \leq \varepsilon r^n,
\end{equation}
is Carleson-prevalent.
\end{definition} 

In Definition \ref{d3}, we were able to require the approximation to go all the way to the boundary,
because, in particular, we had the uniform convergence of Green functions to their limit on compacta of $\RR^n$ rather than just compacta in $\Omega_\infty$,
but this is no longer the case with the gradients, and we shall content ourselves
with estimates on the values of $\nabla G$ on Whitney cubes. We expect better 
estimates than the ones we will prove to be valid, but we shall not pursue this issue here
because the point of this paper is rather to prove weak estimates with compactness arguments.
It will be interesting to investigate more the precise way our approximations quantify, including near the boundary.

\begin{theorem} \label{t41}
Let $\Omega$, $E$, and $L$ satisfy the assumptions of Theorem \ref{t3a1} or
Theorem \ref{t3a1bis}, and let $G$ denote the Green function for $L$ in $\Omega$,
with pole at $\infty$. Then $\nabla G^\infty$ is prevalently close to the gradient of the 
distance to a plane and to
$\nabla D_{\beta,\nu}$, for any  $\beta >0$ and any AR measure $\nu$ on $E$. 
\end{theorem}

\bp We will keep the same structure for the proofs, but will need to replace some estimates. 
We start as in the proof of Theorem \ref{t3a1}, but need to replace 
$\cG_{Gd}(\varepsilon,M)$ with $\cG_{\nabla Gd}(\varepsilon,M)$
in Lemma \ref{l3a2}. This means that instead of \eqref{3a10}, we now want to prove 
\begin{equation} \label{4a6}
\int_{W_M(x,r)} |\nabla\dist(X,P) - \nabla G_k(X)|^2 \leq \varepsilon r^n. 
\end{equation}
In the proof of Theorem \ref{t3a1}, 
the desired estimate
comes from the uniform convergence in Theorem \ref{l2a4}; in the present case, 
 Lemma \ref{l2a2} provides exactly the $L^2$ convergence estimate
that we need (because $G_\infty(X) = \dist(X,P)$ here).

Then we stay in Case 1, as in Theorem \ref{t3a1}, but approximate $\nabla G$
with $\nabla D_{\beta,\mu}$. We start the same way, but 
need to replace our main estimates for the difference between 
$D_{\beta,\mu}(X)$ by $C \dist(X,P)$ by estimates on the difference of gradients.
That is, we need to replace \eqref{3d17}, which itself comes from estimates
for differences of integrals, which culminate with \eqref{3c38}, at least
in the region $H$ of \eqref{3c30}. Here we don't care about what happens
on $B(r,Nr/5) \sm H$ which was treated at the end of the argument, because
this region does not meet the given Whitney region $W_M(x,r)$
if $N$ is large enough and $\eta_2$ is small enough. That is, if we only care about 
$W_M(x,r)$, \eqref{3c38} is enough for the corresponding analogue of \eqref{3d17}.

Now Remark \ref{r31} says that instead of \eqref{3c38} we can prove
\eqref{3c42}, which is just similar but controls the difference of gradients.
That is, repeating \eqref{3c42}, we claim that we have the better estimate
\begin{equation} \label{4a7} 
\begin{aligned}
|\nabla D_{\beta,\mu}(X) - \nabla D(X)|  
\leq C \dist(X,P)^{2+\beta} N^{-\beta-1} r^{-\beta-2},
\end{aligned}
\end{equation}
where $D(X) = (ca_\beta)^{-1/\beta}\dist(X,P)$ as needed (see the line below \eqref{3c34}),
and which is valid in the region
\begin{equation} \label{4a8}
H(N,\eta'_2) = \big\{ X\in B(x,Nr/5) \, ; \, \dist(X,P) \geq C \eta'_2 r \big\}
\end{equation}
where $\eta'_2 = C_1 N \eta^{1/(d+3+\beta)}$ and $\eta_1 = CN \eta^{1/(d+1)}$. 
This estimate will be a good replacement for \eqref{3d17}, but let us check that we do 
not go wrong with the management of Whitney regions.

Recall that we want to prove that $\nabla G$ is prevalently close to 
$\nabla D_{\beta,\mu}$. We give ourselves $\varepsilon > 0$, $M \geq 1$, 
and we need to show that for $(x,r) \in E \times (0,+\infty)$, 
$(x,r) \in \cG_{\nabla GD_\beta}(\varepsilon, M)$ 
prevalently, which means that we want to prove \eqref{4a4}.
We may assume that $(x,r) \in \cG_{\nabla Gd}(\varepsilon/2,M)$, because this is 
a prevalent condition. So we have \eqref{4a3}, and we only need to check that 
\begin{equation} \label{4a9}
\int_{W_M(x,r)} |\nabla \dist(X,P) - c \nabla D_\beta(X)|^2 \leq \varepsilon r^{n}/2.
\end{equation}
We may also assume that $(x,r) \in \cG_{\alpha}(\eta,N)$,
with values of $\eta$ and $N$ that we can choose, and then we get \eqref{4a7},
but maybe for some other $d$-plane $P'$ that comes from the $\alpha$-number.
Let us first assume that $P'=P$, and check that \eqref{4a7} then implies \eqref{4a9}.
In particular we need to make sure that in \eqref{4a7},
$C \dist(X,P)^{2+\beta} N^{-\beta-1} r^{-\beta-1} \leq \varepsilon/2$
when $X \in B(x,Mr)$; we can do this by taking $N$ large enough, depending on $M$ 
and $\varepsilon$. We also need $H(N,\eta'_2)$ above to contain $W_M(x,r)$, 
but this is easy to arrange by taking $N$ large and then $\eta$ small. 
So \eqref{4a7} would imply \eqref{4a9} if we had the same plane.

Maybe $P$, that we get for our proof that $(x,r) \in \cG_{\nabla Gd}(\varepsilon/2,M)$
prevalently, and $P'$, that we choose in terms of $\alpha$-numbers 
(because $(x,r) \in \cG_{\alpha}(\eta, N)$), are not the same. But we can choose them
as close to each other as we want, in particular because the $\alpha$-numbers also control
the flatness, by \eqref{3a17} and \eqref{3a17bis}. See an analogous argument in the proof of Theorem \ref{t3a1}. Then we can control the difference
between $\nabla \dist(X,P)$ and $\nabla D(X)$, where $D$ is the distance associated to
$\dist(X,P')$, and this is enough for \eqref{4a9}.
This ends the proof for Case 1.

\ms
Now we need to follow the proof of Theorem \ref{t3a1bis}, i.e., take care of Case 2.
In this case, we can follow the proof above and, when it comes to applying 
Theorem \ref{l2a4}, we observe as in Case 1 that Lemma~\ref{l2a2} also gives
the $L^2$ convergence of the $\nabla G_k$ to $G^\infty$ (a distance function)
away from $E_\infty = P$. Finally the approximation by $\nabla D_{\beta,\nu}$ is treated 
exactly as in Case 1 (we made sure not to use anything specific, just the good approximation of 
$\nabla D_{\beta,\nu}(X)$ by $\nabla \dist(X,P)$.
\qed

\begin{remark} \label{r41} Assume that $A_k=I$ in Case 2 
(this is the center of our interests anyway). Then
the coefficients of our operator $L_{\alpha,\mu}$ are smooth away from $E$, 
so we could replace the $L^2$
norm in \eqref{4a3} and \eqref{4a4} by a $L^\infty$ norm,
and even write prevalent approximation properties with higher gradients. 
We decided not to check any details here, and anyway the best way to prove such estimates 
is probably to start from the fact that $G$ or its first gradient is close to $D(X) = \dist(X,P)$,
or $\nabla D$, and then follow the same route as for Lemma~\ref{l2a2}. 
That is, we would say that in balls $B(Y,R) \subset \Omega_\infty$, $G$ is prevalently close
to (some) $D$, that $L_{\alpha,\mu} G = 0$ while $\Delta D = 0$, that these two
operators are prevalently close to each other (including for the derivatives of coefficients),
and use some integrations by parts to conclude that  
the derivatives of $G$ are also
close to the derivatives of $D$ in, say, $B(Y,R/2)$. 
The difference between the derivatives of $D_{\beta,\nu}$ and those of $D$ can be 
handled by Remark~\ref{r31}.
\end{remark}

\section{Local variants with $G^X$}
\label{S5}

Maybe the reader does not like too much the function $G^\infty$,
and prefers to use Green functions $G^Y$ with poles at finite distance.
In this section we prove local estimates with the $G^Y$, which will
follow from the previous proofs because the difference between
$G^Y$ and $G^\infty$ is quite small, in particular in small balls
$B(x,r) \subset B(x_0,r_0)$, with $r << r_0$.

We need a local notion of prevalence.

\begin{definition} \label{d51}
Let $E \subset \R^n$ be Ahlfors regular of dimension $d$,
and let $B_0 = B(x_0,r_0)$ be a ball centered on $E$.
We say that the set
$\cG \subset (E \cap B_0) \times (0,r_0)$ is \ub{locally prevalent in $B_0$}
when its complement relative to $B_0$,  $(E \cap B_0) \times (0,r_0) \sm \cG$,
satisfies the Carleson Packing condition \eqref{1a6}. Then 
the (best) constant 
$C$ in \eqref{1a6} is also called the prevalence constant for $\cG$ in $B_0$.
\end{definition}

It may happen that $\cG$ comes from a globally defined subset of $E \times (0,+\infty)$,
but this will not always be the case, typically because $\cG$ depends on the choice of
a pole that we take far enough from $E \cap B_0$, but anyway when this is the case we simply
forget about that part of $\cG$ that does not lie in $(E \cap B_0) \times (0,r_0)$.

More importantly, we usually ask for uniform estimates, by which we mean that the Carleson constant $C$ in \eqref{1a6} does not depend on $B_0$ (and the associated set $\cG$).

Now what sets $\cG$ (associated to $B_0$) shall we consider? 
Let $C_0 \geq 1$ be given; we shall only consider poles $Y$ such that
\begin{equation} \label{5a2}
\dist(Y, E \cap B_0) \geq C_0^{-1} r_0
\end{equation}
and then, for the local prevalent approximation of $G^Y$ by a distance function, we 
select a pole $Y$ such that \eqref{5a2} holds, and use the good set 
$\cG_{Gd}^{B_0,Y}(\varepsilon, M)$ of pairs $(x,r) \in (E \cap B_0) \times (0,r_0)$ 
such that there exists a $d$-plane $P(x,r)$ and a positive constant $c > 0$, with
\begin{equation} \label{5a3}
|\dist(X,P) - c G^Y(X)| \leq \varepsilon r
\text{ for } X \in \Omega \cap B(x,Mr).
\end{equation}

For the approximation of $G^Y$ by 
a distance function $D_{\beta,\nu}$, 
we use instead the good set $\cG_{GD_\beta}^{B_0,Y}(\varepsilon, M)$ defined the same way, but with 
\eqref{5a3} replaced by 
\begin{equation} \label{5a4}
|D_{\beta,\nu} 
- c G^Y(X)| \leq \varepsilon r
\text{ for } X \in \Omega \cap B(x,Mr).
\end{equation}
Finally we define the prevalent approximation property in the way that will suit our
results.

\begin{definition} \label{d52}
Let $E \subset \R^n$ and $L$ be as above. We say that the
$G^Y$ are locally prevalently close to the distance to a plane
(respectively, locally prevalently close to the function 
$D_{\beta,\nu}$)
when for each ball $B_0 = B(x_0, r_0)$ and each $Y \in \Omega$ such that \eqref{5a2}
holds, the set $\cG_{Gd}^{B_0,Y}(\varepsilon, M)$
(respectively, $\cG_{GD_\beta}^{B_0,Y}(\varepsilon, M)$)
is locally prevalent in $B_0$, with constants that do not depend on $B_0$ or 
on $Y$ satisfying \eqref{5a2}.
\end{definition}

Notice that this definition depends on $C_0$; this will not disturb us, because 
we will be able to prove the desired result for any $C_0$, and with Carleson constants
that do not even depend too badly on $C_0$. In fact, with a little bit of manipulations
of the definitions, it is easy to check that different constants $C_0$
give the same notion of local prevalence (but different constants). 
The proof below gives a good idea of that.

\begin{theorem} \label{t51}
Retain the assumptions of Theorem~\ref{t3a1} or \ref{t3a1bis}. 
Denote by $G^Y$, $Y \in \Omega$, the Green function associated to $L$
in $\Omega$ with pole at $Y$.
Then the $G^Y$ are locally prevalently close to the distance to a plane,
and to any function $D_{\beta,\nu}$ where $\beta > 0$ and $\nu$
is an AR measure on $E$.
\end{theorem}

As before, in the second statement, we are even allowed to define $D_{\beta, \nu}$
with a different AR measure $\nu$ than the one in the definition of $L_\alpha$.

The reader should not be too surprised that we also allow poles $Y$ that are far
from $E \cap B_0$ but close to (faraway parts of) $E$; the comparison principle says
that the corresponding Green functions do not behave in a different way than if we forced
$Y_k$ to be corkscrew points in large balls centered at $x_0$. 
Of course we may need to normalize differently,
i.e., choose very different constants $c$ in \eqref{5a3} or \eqref{5a4}, but this is all right.

We do not expect to have a good approximation of $G^Y$ in the balls $B(x,r)$ such that
$x\in E \cap 2B_0$ and $r \geq C^{-1} r_0$, even when $E$ is a plane and $L = \Delta$,
but those balls will satisfy a Carleson packing condition. Indeed,
once $\varepsilon$ and $M$ are given
(as in \eqref{5a3} or \eqref{5a4}), our first action is to remove the set
\begin{equation} \label{5a7}
\cB_0 = \big\{ (x,r) \, ; \, x \in E \cap B_0 \text{ and } \tau r_0 \leq r \leq r_0 \big\}
\end{equation}
for some very small $\tau \in (0,1)$ that we allow to depend on 
$\varepsilon$, $M$, and $C_0$. 
We can do this because $\cB_0$ satisfies a Carleson packing condition
(with a constant that depends on $\tau$ and the AR constants).

Now we pick $(x,r) \in (E \cap B_0) \times (0, \tau r_0)$, and follow the same argument as
in the proofs above. Essentially, we use the fact that $r$ is so small that $Y$ will appear as very far, and seen from
$B(x,r)$, $G^Y$ will look a lot as $G^\infty$. Let us for instance consider Case 1, for the approximation by distances
to planes. There is a moment, in the proof of Lemma~\ref{l3a2},
where we want to show by compactness that if $(x,r)$ satisfies some prevalently good properties,
then $(x,r) \in \cG_{Gd}^{B_0,Y}(\varepsilon, M)$. We approach it by contradiction and
find a sequence of counterexamples $(\Omega_k, E_k, L_k, x_k, y_k)$, and now
we also need to add a larger ball $B_{0,k} = B(x_{0,k}, r_{0,k})$ and a pole
$Y_k$ such that $\dist(Y_k, E_k \cap B_{0,k}) \geq C_0^{-1} r_{0,k}$.
We also let $\tau_k$ tend to $0$, and recall that $r_k \leq \tau_k r_{0,k}$.
Finally, $G_k$ is no longer a Green function at $\infty$, but rather
$G_k = G_{L_k}^{Y_k}$, with the pole $Y_k$.

As before we normalize things so that $x_k=0$ and $r_k = 1$.
Then $|Y_k| \geq C^{-1} r_{0,k} \geq C_0^{-1} \tau_k^{-1}$ tends to $+\infty$,
and by Theorem \ref{l2a4} and Remark \ref{r2a18}, the functions $c_k G_k^{Y_k}$
(correctly normalized by constants $c_k$) converge to the Green function 
$G_\infty^\infty$ for $\Omega_\infty$ and $L_\infty$. From there on,
we can keep the same proof. Case 2 is similar, and for the approximation
with our special distances $D_{\beta,\nu}$, we only need to observe that since
$D_{\beta,\nu}$ and the distance to a plane are prevalently close to each other, the result
stays true locally, just by definitions.
We leave the details to the reader.
\qed

\begin{theorem} \label{r52}
Retain the assumptions of Theorem~\ref{t3a1} or \ref{t3a1bis}. 
Denote by $G^Y$, $Y \in \Omega$, the Green function associated to $L$ in $\Omega$ 
with pole at $Y$.
Then the $\nabla G^Y$ are locally prevalently close to the gradient of the distance to a plane,
and to any $\nabla D_{\beta,\nu}$ where $\beta > 0$ and $\nu$ is an AR measure on $E$.
\end{theorem}

The proof would follow the proof of Theorem \ref{t41}, using the fact that 
 Lemma \ref{l2a2} remains valid in the conditions of Remark \ref{r2a18}. We could also
deduce the approximation of the gradients from the approximation of the functions,
as in the proof of Lemma \ref{l2a2} itself.

\section{The basic converse, with the distance to a $d$-plane}
\label{S6}

We want to say than when $\Omega \subset \R^n$
is bounded by an unbounded AR set $E = \d\Omega$ of integer dimension $d$,
and the Green function $G^\infty$ associated to one of our operators
is prevalently close to the distance to a $d$-plane,
then $E$ is uniformly rectifiable.

Let us 
state the assumptions of our theorem in advance.
We are given a domain $\Omega \subset \R^n$, whose boundary $E = \d \Omega$
is an (unbounded) AR set of integer dimension $d \leq n-1$.
When $d=n-1$, we also assume that $\Omega$ contains corkscrew balls and 
Harnack chains, as in \eqref{2b4}. Recall that when $d < n-1$, we don't need 
to ask for \eqref{2b4}, because it is always true.

Then we are also given an operator $L = - \dv A \nabla$, and we distinguish two cases.
In Case~1, we assume that $d = n-1$ and that $A$ satisfies the ellipticity conditions 
\eqref{2b2} and \eqref{2b3}. In particular we do not need to know that
$L$ is close to a constant coefficient operator in any way:  as in the next case,
we just want to make sure that there is a reasonable definition of the Green function
$G^\infty$, as in Section \ref{S2}, with some basic properties.

In Case 2, we allow any integer $0 < d \leq n-1$, and any operator \eqref{2b6-bis}, 
or, equivalently, $L = -\dv A\dist(X,E)^{d-n+1} \nabla$,  such that $A$ satisfies the ellipticity
conditions \eqref{2b2} and \eqref{2b3}. 
This includes the good operators $L_\alpha$ of \eqref{1a2}-\eqref{1a3}, 
but also any degenerate elliptic operator in the class that was studied in 
\cite{DFM2}. Again we need something like this because 
we want to quote general results concerning Harnack's inequality and the 
Green function, but no delicate geometry will be used concerning $L$.

Since $d$ is integer in the present Theorem, Case 2 automatically includes Case 1.

\begin{theorem} \label{t61}
Let $\Omega \subset \R^n$ be a domain bounded by an unbounded $d$-AR set $E=\po$, $0<d<n$ integer, and assume that when $d= n-1$ the domain $\Omega$ contains corkscrew balls and  Harnack chains, as in \eqref{2b4}. Let $L = -\dv A\dist(X,E)^{d-n+1} \nabla$,  with $A$ elliptic as per \eqref{2b2} and \eqref{2b3}.
Assume that the 
Green function $G^\infty$ associated to $\Omega$ and $L$,
with a pole at $\infty$, is prevalently close to the distance to a $d$-plane. 
Then $E$ is uniformly rectifiable.
\end{theorem}

\bp Our assumption is that for every choice of $\varepsilon > 0$ and $M \geq 1$,
the good set $\cG_{Gd}(\varepsilon, M)$ of Definition \ref{d3} is prevalent.
We want to show that for each $\tau > 0$, the other good set $\cG_{ur}(\tau)$ of
Definition \ref{d2} is prevalent, and the simplest is to show that if $M$ and $\varepsilon$
are chosen correctly, $\cG_{Gd}(\varepsilon, M) \subset \cG_{ur}(\tau)$.
In the present case, taking $M=10$ will be enough.

So we pick $(x,r) \in \cG_{Gd}(\varepsilon, 10)$ and by definition there is a $d$-plane $P$
and a constant $c>0$ such that, as in \eqref{1a10}, 
\begin{equation} \label{6a2}
|\dist(X,P) - c G^\infty(X)| \leq \varepsilon r
\text{ for } X \in \Omega \cap B(x,10r). 
\end{equation}
Recall that $G^\infty$ is H\"older continuous on $\Omega \cap B(x,10r)$, so
we can safely say that $\lim_{X \to y} G^\infty(X) = G^\infty(y) = 0$
when $y \in E \cap B(x,10r)$ without thinking about 
the notion of trace, and by \eqref{6a2}
\begin{equation} \label{6a3}
\dist(y,P) \leq \varepsilon r \text{ for } y \in E \cap B(x,10r).
\end{equation}

If we have this property alone, we get what is known as the weak geometric lemma,
which is somewhat weaker than the uniform rectifiability of $E$, so we have to continue 
the argument and show that conversely, 
\begin{equation} \label{6a4}
\dist(X,E) \leq \tau r/2 \text{ for } X \in P \cap B(x,r).
\end{equation}
Obviously, if we prove \eqref{6a4} and take $\varepsilon < \tau/2$, we get that 
$d_{x,r}(E,P) \leq \tau$, and hence $(x,r) \in \cG_{ur}(\tau)$ as needed.

Set $u = c G^\infty$. We will derive \eqref{6a4} via 
estimates on the size of $u$. Since \eqref{6a3} says that 
$\dist(y,P) \leq \varepsilon r$ for $y \in E \cap B(x,10r)$, 
we now get that
\begin{equation} \label{6a5}
u(X) \leq \dist(X,P) + \varepsilon r \leq 11 r \text{ for } X \in \Omega \cap B(x,10r),
\end{equation}
and we shall now find a point $A_1$ such that $u(A_1)$ is not too small.
First take a corkscrew point $A_0$ for $B(x,r)$. That is, $A_0 \in B(x,r) \cap \Omega$
and $\dist(A_0,E) \geq C^{-1} r$, and then \eqref{6a3} allows us to find a new point $A_1$
such that $|A_1-A_0| \leq (2C)^{-1} r$ (and so $A_1\in \Omega \cap B(x,2r)$)
and in addition $\dist(A_1,P) \geq (2C)^{-1} r$. 
Thus we can also use $A_1$ as a corkscrew point, and now 
\begin{equation} \label{6a6}
u(A_1) \geq (2C)^{-1} r - \varepsilon r \geq (3C)^{-1} r
\end{equation}
by \eqref{6a2} (and if $\varepsilon$ is chosen small enough). Now we claim that
\begin{equation} \label{6a7}
u(X) \geq C [r^{-1}\dist(X,E)]^\gamma  \, u(A_1)   \geq c' r [r^{-1}\dist(X,E)]^\gamma 
\ \text{ for } X \in \Omega \cap B(x,r),
\end{equation}
where $C \geq 1$, $\gamma > 0$, and then $c'>0$ depend only on the geometric 
constants for $\Omega$ and the ellipticity (or degenerate ellipticity) constant for $L$. 
This will be 
basically due to the fact that the length of the Harnack chain from 
$A_1$ to $X$ is $C \,\log(\dist(X,E)^{-1}r)$.

Indeed let $X \in \Omega \cap B(x,r)$, and pick $y\in E \cap B(x,2r)$ such that
$|X-y| = \dist(X,E)$. Then pick for $k \geq 0$ a corkscrew point $X_k$ for 
$B(y,2^k\dist(X,E))$. We can take $X_0 = X$, and we stop as soon as $2^k\dist(X,E) \geq  r$.
Thus our final $X_k$ lies in $B(y,2r) \subset B(x,4r)$, and its distance to $E$ is more than
$C^{-1} r$. By Harnack's inequality, $u(A_1) \leq C u(X_k)$. We can also
apply Harnack's inequality to consecutive points of the sequence, and get that
$u(X_{j+1}) \leq C u(X_j)$ for $0 \leq j < k$. 
Altogether, $u(A_1) \leq C^{k+1} u(X_0) = C^{k+1} u(X)$; \eqref{6a7} follows because
$k+1 \sim \log_2(\dist(X,E)^{-1}r)$.

Now set 
\begin{equation} \label{6a8}
H = \big\{ X \in B(x,2r) \, ; \, \dist(X,P) \leq 2\varepsilon r \big\}
\end{equation}
and let $X \in H \cap \Omega$ be given. By \eqref{6a2},
$$
u(X) =c G^\infty(X) \leq \dist(X,P) + \varepsilon r \leq 3\varepsilon r\quad \mbox{for all} \quad X\in H.
$$
Then by \eqref{6a7}, 
$$c' r [r^{-1}\dist(X,E)]^\gamma \leq u(X)\leq  3\varepsilon r,$$ 
hence
$r^{-1}\dist(X,E) \leq C \varepsilon^{1/\gamma}$; we choose $\varepsilon$ so small that
$C \varepsilon^{1/\gamma} \leq \tau/3$, and get that 
\begin{equation} \label{6a9}
\dist(X,E) \leq \tau r/3
\ \text{ for } X \in H \cap \Omega \cap B(x,2r).
\end{equation}
When $d < n-1$, $E$ is so small that $\R^n \sm E$ is connected,
so $\R^n = \Omega \cup E$. Then for any $X \in P \cap B(x,2r)$, either
$X \in \Omega$ and $\dist(X,E) \leq \tau r/3$ by \eqref{6a9}, or else 
$X \in E$ and $\dist(X,E) =0$. So \eqref{6a4} holds in this case.

We are left with the case when $d=n-1$. Recall that $\dist(A_1,P) \geq (2C)^{-1} r$,
so $A_1$ lies in $B(x,2r) \sm H$. Denote by $V$ the component of $B(x,2r) \sm H$
that contains $A_1$; then $V \subset \Omega$ because it meets $\Omega$
and not $E$ (by \eqref{6a3} and \eqref{6a8}).

Now let $X \in P \cap B(x,r)$ be given. There exists $Y \in V$ such that 
$|Y-X| \leq 3\varepsilon$, and by \eqref{6a9} $\dist(Y,E) \leq \tau r/3$;
if $\varepsilon$ is small enough, this implies that $\dist(X,E) \leq \tau r/2$;
this proves \eqref{6a4} and the Theorem \ref{t61} follows.
\qed

\begin{remark} \label{r61}
Notice that in the course of the proof of Theorem \ref{t61} we hardly use the fact that
$G^\infty$ is a Green function (we only use the fact that it vanishes on $E$ and the nondegeneracy property \eqref{6a7}). This, however, is not shocking -- see, e.g., \cite{DEM}, where the free boundary results also rely very moderately on the particular properties of the distance function. The situation may be even more peculiar in the special case $d=n-1$ 
because it could happen (but the authors did not manage to
prove this) that if $E$ is Ahlfors regular 
of dimension $n-1$ and satisfies the Weak Geometric Lemma (see for instance \cite{DS}), and, in addition, $\Omega$ contains corkscrew balls and 
Harnack chains, then $E$ is uniformly rectifiable. If this were true, then we could stop
the argument as soon as we proved \eqref{6a3}, and the only properties of $G^\infty$
that we would 
use is the prevalent good approximation and the fact that $G^\infty = 0$ on $E$.
\end{remark}

\begin{remark} \label{r62}
Theorem \ref{t61} also holds when we assume that the Green functions $G^Y$
are locally prevalently close to distances to a $d$-plane (see Definition \ref{d52}).
The proof is the same; that is, we only used the fact that the function $u = c G^\infty$
is $L$-harmonic in $\Omega \cap B(x,10r)$, and reducing to the case when
the pole $Y$ lies outside of $B(x,10r)$ is just a matter of skipping the large balls
$B(x,r)$, $r \geq r_0/C$, i.e., restricting our attention to
pairs $(x,r) \in (E \cap B(x_0,r_0)) \times (0, r_0/C)$, in the Definition \ref{d51}
of local prevalence.
\end{remark}

\begin{remark} \label{r63}
Theorem \ref{t61} also holds when we replace $\dist(X,E)$ with any power $\dist(X,E)^\alpha$,
with the same proof. This is a quite natural condition 
as the Green function is expected to be only H\"older continuous at the boundary, that is, to behave as a power of the distance.
\end{remark}

\begin{remark} \label{r64}
In this section and the next ones, we give statements with distance functions, 
and not their gradients. The results of these sections should also hold with 
the prevalent approximation with the gradients, but we decided not to check them.
Note however that in principle conditions on the good approximation of $\nabla G$ are 
stronger than those on $G$ itself, because $G$ can essentially be computed from its gradient.
\end{remark}

\section{The converse with the distance $D_\alpha$.}
\label{S7}

We are now interested in finding out whether if the Green function $G^\infty$ associated 
to one of our operators is prevalently close to some $D_\alpha$ (and later in this section,
to some power of $D_\alpha$), then $d$ is an integer and $E$ is uniformly rectifiable.

The situation is more challenging than with the distance to $d$-planes, first
because we allow dimensions that are not integers (this made little sense when discussing the distance to a plane), and also,
due to the fact that $D_\alpha$ vanishes on $E$, and not on a $d$-plane
as in Section \ref{S6}, 
we cannot simply use the fact that $D_\alpha$ is small near $E$; that is,
the answer is no longer  hidden in the question.

Contrary to the previous section, we are now restricted to our special classes of elliptic operators, satisfying  \eqref{3b3} or built from the distance functions. 
We start with the most classical case of co-dimension $1$.

\begin{theorem} \label{t71}
Let $\Omega \subset \R^n$ be such that $E = \d \Omega$ is Ahlfors regular
of dimension $n-1$ and assume, in addition, that the domain $\Omega$ has interior corkscrew balls and Harnack chains, that is, \eqref{2b4} holds. Let $L = - \dv A \nabla$
satisfy \eqref{2b2}, \eqref{2b3}, and \eqref{3b3}.
Suppose that the Green function $G^\infty$ for $L$ on $\Omega$
is prevalently close to $D_{\alpha,\mu}$ for some $\alpha > 0$ and some AR 
measure $\mu$ of dimension $n-1$ on $E$. Then $E$ is uniformly rectifiable.
\end{theorem}

\bp This time we shall not use Definition \ref{d2} to prove the uniform rectifiability, 
but rather show that, in addition to being 1-sided NTA, $\Omega$ has exterior corkscrew points. 
This is also known as \ub{Condition B}. We say that the (unbounded) AR set $E$
of dimension $n-1$ satisfies Condition B when there is a constant $c > 0$ such that
for $x\in E$ and $r > 0$, 
\begin{equation} \label{7a8}
\begin{aligned}
&\text{we can find $X_1, X_2 \in B(x,r)$, that lie in different connected components } 
\cr&\hskip1.5cm
\text{of $\R^n \sm E$, and such that $\dist(X_i,E) \geq c r$ for $i=1,2$.}
\end{aligned}
\end{equation}
This condition was introduced by S. Semmes in \cite{Se1}, and is known to imply the uniform
rectifiability of $E$ (see \cite{Se1} in the smooth case and, probably for the simplest proof,
\cite{DJ}).

This property can be also written in terms of prevalent sets. Let
$\cG_{CB}(c)$ denote the set of pairs $(x,r) \in E \times (0,+\infty)$ such that \eqref{7a8} holds. The definition requires that $\cG_{CB}(c) = E \times (0,+\infty)$, but we claim that 
if the set $\cG_{CB}(c)$ is prevalent for some $c > 0$, then $E$ satisfies Condition B 
(with a worse constant). 

Indeed, if $\cG_{CB}(c)$ is prevalent, we first claim that there is a constant $a > 0$
such that for each $(x,r) \in E \times (0,+\infty)$, we can find
$y\in E \cap B(x,r/2)$ and $t \in (a r, r/2)$ such that $(y,t) \in \cG_{CB}(c)$.
Otherwise the set $E \cap B(x,r/2) \times (a r, r/2)$ is contained in the
bad set $\cB = E \times (0,+\infty) \sm \cG_{CB}(c)$, and then
\begin{eqnarray} \label{7a9}
\int_{y \in E \cap B(x,r)}\int_{0<t<r} \1_{\cB}(y,t) \frac{d\mu(y) dt}{t} 
&\geq& \int_{y \in E \cap B(x,r/2)}\int_{ar<t<r/2}  \frac{d\mu(y) dt}{t} 
\nn\\
&\geq& \mu(E \cap B(x,r/2) \int_{ar<t<r/2} \frac{dt}{t} \geq C^{-1} r^{d} \ln(1/(2a)),
\end{eqnarray}
which contradicts \eqref{1a6} if $a$ is chosen small enough, depending on
the Carleson constant for $\cB$. But once we obtain that pair $(y,t)$, the two points
$X_1$ and $X_2$ that we get from \eqref{7a8} for $B(y,t)$ also work for
$B(x,r)$, although with the worse constant $ac$.

Return to the theorem, and let $\Omega$, $E$, $L$ satisfy the assumptions.
We only need to find $c > 0$ such that $\cG_{CB}(c)$ is prevalent.
Since we know that for each choice of $\varepsilon$ and $M$, 
the sets $\cG_{GD_\alpha}(\varepsilon, M)$ and $\cG_{cc}(\varepsilon, M)$
are prevalent, it is enough it show that 
\begin{equation} \label{7a7}
\cG_{GD_\alpha}(\varepsilon, M) \cap \cG_{cc}(\varepsilon, M) \subset \cG_{CB}(c)
\end{equation}
for some choice of $c>0$, $M$, and $\varepsilon$.

We proceed by contradiction, assume that \eqref{7a7} fails for
$c_k = \varepsilon_k = M_k^{-1}= 2^{-k}$, and pick a counterexample for each $k$.
This means, an open set $\Omega_k$ bounded by an AR 
set $E_k$, an operator $L_k$, and a pair $(x_k, r_k) \in E \times (0,+\infty)$
that all satisfy the assumptions, and for which
$(x_k, r_k) \in \cG_{GD_\alpha}(\varepsilon_k, M_k) \cap \cG_{cc}(\varepsilon_k, M_k) 
\sm \cG_{CB}(c_k)$.

By translation and dilation invariance, we may assume that $x_k = 0$ and $r_k = 1$. 
We may also extract a subsequence so that $\Omega_k$
converges to a limit $\Omega_\infty$, $E_k$ converges to $E_\infty$,
and per Theorem \ref{l2a4}, $\Omega_\infty$ is bounded
by $E_\infty$ and $\Omega_\infty$ satisfies the assumption \eqref{2b4}. Moreover,
we can extract a further subsequence so that the AR measure $\mu_k$
given on $E_k$ converges weakly to a measure $\mu$, and this limit is
automatically AR. Because $(0,1) \in \cG_{cc}(2^{-k}, 2^k)$,
the matrix $A_k$ of $L_k$ converges, in $L^1_{loc}(\Omega_\infty)$, 
to a constant matrix $A_0$.

By Theorem \ref{l2a4}, the Green function $G_k = G_k^\infty$
for $L_k$ in $\Omega_k$, correctly normalized, converges to
the Green function $G = G_\infty^\infty$ for $L_0$ on $\Omega_\infty$,
uniformly on compact sets of $\Omega_\infty$. Here $L_0$ is the
constant coefficients operator $L_0 = - \dv A_0 \nabla$.
Hence, using the facts that $(0,1) \in \cG_{GD_\alpha}(2^{-k}, 2^k)$
and  $\mu_k$ stays uniformly AR and converges to $\mu$,
we can conclude that $G = C D_{\alpha, \mu}$. 

We shall start many of our ``free boundary" compactness arguments
in this fashion, but now we need to use the specific assumption of co-dimension $1$. 
Why does 
 the fact that $G=CD_{\alpha, \mu}$ imply 
 something nice about the domain?
The simplest path would be to brutally use a known fact that on 1-sided NTA domains 
if  the harmonic measure (associated to a constant coefficient elliptic operator on $\Omega$) is absolutely continuous 
with respect to $\mu$ (any AR measure on $E$) and given by an $A_\infty$ weight, 
then $E$ is uniformly rectifiable. This is essentially due to \cite{HMU}. A careful reader might notice that formally speaking, they treat the Laplacian only, but this extends to all constant coefficients symmetric elliptic operators, and then a symmetrization argument could be used to get rid of the assumption of symmetry. See, e.g., \cite{HMMTZ} where an analogous result is obtained in much bigger generality (so it could serve as a reference by itself), but in particular, the symmetrization is discussed in the beginning of Proof of Theorem~1.6.

Our limiting domain $\Omega_\infty$, with the constant coefficients operator $L_0$,
admits the Green function $G =C D_{\alpha, \mu}$, and in particular
\begin{equation} \label{7b5} 
C^{-1} \dist(X,E) \leq G(X) \leq C \dist(X,E)
\ \text{for } X \in \Omega_\infty.
\end{equation}
We claim that this implies that the harmonic measure $\omega^\infty$, 
with a pole at $\infty$, is also comparable to $\H^{n-1}_{\vert E_\infty}$, 
in the sense that 
$$
C^{-1}\H^{n-1}_{\vert E_\infty} \leq \omega^\infty 
\leq C \H^{n-1}_{\vert E_\infty}\quad \mbox{on} \quad E_\infty.
$$
Once this claim is established,  we can apply the aforementioned result and obtain that 
\begin{equation} \label{7b6}
\text{$E_\infty$ is uniformly rectifiable.}
\end{equation}

The claim is due to the following estimate for the Green function, valid on any 1-sided NTA domain with an unbounded $d$-dimensional Ahlfors regular boundary when $d\geq n-1$, and on any domain with an unbounded $d$-dimensional Ahlfors regular boundary when $d< n-1$.
We consider the Green function with pole at a faraway point $Y$
for any elliptic operator $L$ on $E$, with \eqref{2b2} and \eqref{2b3}, 
or the operator given by \eqref{2b6}.
Under these assumptions, there is a constant $C \geq 1$,
that depends on the dimensions and the constants in our assumptions on $\Omega$ and $L$, 
such that if $B = B(x,r)$ is a ball centered on $E = \d\Omega$, and 
$A^B \in B \cap \Omega$ denotes a corkscrew point for $B$, 
\begin{equation} \label{7b7}
C^{-1} r^{n-2} G^{Y}(A^B) \leq \omega^Y(E \cap B) \leq C r^{n-2} G^{Y}(A^B)
\end{equation}
as soon as the pole $Y$ lies out of $2B$. We refer to Lemma 15.28 in \cite{DFM5},
but in co-dimension 1 this is of course a classical result. 
To be careful with the normalization of the Green functions, let us 
consider balls $B$ that are contained in $B_0$ for some fixed ball $B_0$ and pick any 
$Y \in \Omega \sm 2B_0$. 
The comparison principle says that  $C^{-1} G^{Y} \leq G^\infty \leq C G^Y$, 
with a uniform constant, if for instance we normalize
$G^\infty$ by $G^\infty(A_0) = G^Y(A_0)$ for a corkscrew point $A_0$ for $B_0$.
Then by  \eqref{7b5}
\begin{equation} \label{7b8}
C^{-1} r \leq  c\,G^{Y}(A^B) \leq C r,
\end{equation}
where the normalizing constant $c$ depends on $B_0$, but not $C$.
By \eqref{7b7}, this yields
\begin{equation} \label{7b9} 
C^{-1} r^{n-1} \leq c \,\omega^Y(E \cap B) \leq C r^{n-1}
\ \text{ for } B \subset B_0,
\end{equation}
which implies by an easy covering argument that $\omega^Y$ is essentially proportional
to $\H^{n-1}$ on $E \cap \frac12 B_0$. So we may apply the result of \cite{HMU},
and we get \eqref{7b6}.

Recall that we want to use \eqref{7b6} to derive a contradiction, in fact with our assumption
that $(0,1) = (x_k, r_k) \notin \cG_{CB}(c_k)$.

Let $\varepsilon > 0$ small, to be chosen soon.  We know that the set $\cG_{ur}(\varepsilon)$ 
associated to $E_\infty$ is prevalent, so by the same argument as in \eqref{7a9} we can find 
$y \in E_\infty \cap B(0,1/2)$ and $t \in (a, 1/2)$ such that
$d_{y,t}(E_\infty,P) \leq \varepsilon$ for some hyperplane $P$.

Denote by $H_+$ and $H_-$ the two connected components of
$H = \big\{ X \in B(y,a) \, ; \, \dist(X,P) \geq 2a\varepsilon r \big\}$, and also let
$A_0$ be a corkscrew point for $\Omega_\infty$ in $B(0,a)$;
we know that it exists because $0 \in E_\infty$ and $\Omega_\infty$ satisfies
\eqref{2b4}. Then $\dist(A_0,P) \geq C^{-1}a$ because $d_{y,t}(E_\infty,P) \leq \varepsilon$,
and so $A_0 \in H_+$, say. 

For $k$ large, $H_+$ and $H_-$ do not meet $E$, because
$d_{y,t}(E_\infty,P) \leq \varepsilon$ and $E_\infty$ is the limit of $E_k$; then
we know that $H_+ \subset \Omega_k$ (because $A_0 \in \Omega_k$ for $k$ large),
and we have two options. The pleasant one is when $H_-$ is contained in some other component
of $\R^n \sm \Omega_k$, because this contradicts immediately our contradiction
assumption that $(0,1) \notin \cG_{CB}(c_k)$ (relative to $\Omega_k$): the two points
$A_0$ and its symmetric relative to $P$ would fit in \eqref{7a8}.
The other option is that $H_-$ meets $\Omega_k$, hence $H_- \subset \Omega_k$
(because it does not meet $E_k$). This is impossible too, this time because
if $\varepsilon$ is small enough, we contradict the existence of Harnack chains in $\Omega_k$;
see the argument below \eqref{3a6}.

We finally proved \eqref{7a7} by compactness and contradiction, and 
Theorem \ref{t71} follows.
\qed

\begin{remark} \label{r71}
It seems that we should not be forced to use \cite{HMU, HMMTZ}: this is a rather indirect route, and we have a detailed information about the Green function at hand, in particular, we know that it is a multiple of a suitably defined distance to the boundary, not just the fact that it satisfies the bounds from above and below \eqref{7b5}.  Yet, these matters are highly non-trivial, both for the Green function and for the harmonic measure. For instance, it is
known that there are cones $E$, other than hyperplanes,
such that $\omega^\infty$ is proportional to $\H^{n-1}_{\vert E}$. 
We do not know (but did not try to make a computation)
whether the Green function will be a multiple of $D_\alpha$ for these examples.  
It could happen that in the situation of $G = CD_{\mu,\alpha}$, $\mu$ 
really must be a flat measure (instead of being merely uniformly rectifiable, as in \eqref{7b6}. 
We will return to this issue in the next section. 
\end{remark}

We now switch to the question of prevalent approximation by $D_{\alpha,\mu}$ of the Green 
functions for elliptic operators in dimensions $d \in (n-2,n) \sm \{n-1\}$. 
We want to say that this does not happen. That is, the fact that the Green function is 
prevalently close to $D_{\alpha,\mu}$ implies that the dimension is integer, for the operators in Theorem~\ref{t71}, and, more generally, for all classical elliptic operators. 

\begin{proposition} \label{t72}
Let $\Omega \subset \R^n$ be such that $E = \d \Omega$ is Ahlfors regular of dimension 
$d \in (n-2,n) \sm \{n-1\}$. If $d > n-1$, suppose in addition that $\Omega$ contains 
corkscrew balls and Harnack chains, as in \eqref{2b4}. Then let $L = - \dv A \nabla$ be an
elliptic operator (as in \eqref{2b2} and \eqref{2b3}), and denote by $G^\infty$
the Green function for $L$ on $\Omega$. Then $G^\infty$ is not prevalently close to 
$D_{\alpha,\mu}$ for any choice of $\alpha >0$ and an AR 
measure $\mu$ of dimension $d$ on $E$.
\end{proposition}

This was demoted to the rank of proposition, because it will appear soon that the approximation 
fails just because $D_{\alpha,\mu}$ does not have the right homogeneity.

\bp Suppose the Proposition fails, so that we can find an example $(\Omega, E, L)$
for which $G^\infty$ is prevalently close to some $D_{\alpha,\mu}$.
This means in particular that we can find pairs $(x,r) \in \cG_{GD_\alpha}(\varepsilon, M)$,
with arbitrary values of $\varepsilon$ and $M$.

Pick such a pair; by translation and dilation 
invariance we may assume that $(x,r) = (0,1)$, so we 
deduce from \eqref{1a12} that there is a constant $c > 0$ such that
\begin{equation} \label{7b12}
|D_{\alpha,\mu}(X) - c G^\infty(X)| \leq \varepsilon 
\ \text{ for } X \in \Omega \cap B(0,M).
\end{equation}
In particular, recalling that $D_{\alpha,\mu}(X) \sim \dist(X,E)$ and if $\varepsilon$ is small enough, 
we get that 
\begin{equation} \label{7b13}
C^{-1} \dist(X,E) \leq c G^\infty(X) \leq C  \dist(X,E)
\end{equation}
for $X \in \Omega \cap B(0,M)$ such that $ \dist(X,E)) \geq 1$.
Now pick a pole $Y \in \Omega \sm B(0, 10M)$ and a  corkscrew point $X_0$ for $B(0,M)$,
then set 
\begin{equation} \label{7b14}
\lambda = G^Y(X_0)/G^\infty(X_0)
\end{equation}
and apply the comparison principle to the two function $G^\infty$ and $G^Y$, 
say, in $\Omega \cap B(0,2M)$. We get that 
\begin{equation} \label{7b15}
C^{-1} \lambda  \leq  \frac{G^Y(X)}{G^\infty(X)} \leq C \lambda 
\ \text{ for } X \in \Omega \cap B(0,M)
\end{equation}
and, because of \eqref{7b13}, 
\begin{equation} \label{7b16}
C^{-1} \dist(X,E) \leq \frac{c}{\lambda} G^Y(X) \leq C \dist(X,E)
\end{equation}
when in addition $\dist(X,E)) \geq 1$. Now consider a ball 
$B = B(x,r) \subset B(0,M)$, and choose a corkscrew point $A^B$ for $B$.
We restrict to the case when $r \geq C_0$, where $C_0$ is chosen
so that \eqref{7b16} always holds for $A^B$
and yields 
$$G^Y(A^B) \sim \frac{\lambda}{c} \dist(X,E) \sim \frac{\lambda}{c} \, r.$$
Then apply \eqref{7b7} to $B$; this yields
\begin{equation} \label{7b17} 
C^{-1} \, \frac{\lambda}{c} \, r^{n-1} 
\leq \omega^Y(E \cap B(x,r)) \leq C \,\frac{\lambda}{c} \, r^{n-1}.
\end{equation}
We claim that this is incompatible with the Ahlfors regularity of $E$.
We first take $B_0 = B(0,M/2)$ and find that
$\omega^Y(E \cap B_0) \sim \frac{\lambda}{c} M^{n-1}$.

First assume that $d > n-1$; we can find more than  $C^{-1} M^d$
disjoint balls $B$ of radius $C_0$ that are contained in $B_0$, and hence
\eqref{7b17} yields 
$$\frac{\lambda}{c} M^{n-1} \gtrsim \,\omega^Y(E \cap B_0) 
\geq \sum_B \omega^Y(E \cap B_0)
\geq C^{-1} \frac{\lambda}{c} M^d,$$ 
a contradiction if $M$ is large enough and $d > n-1$. Observe that
$\frac{\lambda}{c}$ disappears as it should (we could have replaced $G^\infty$ by a multiple
so that $\frac{\lambda}{c} = 1$ anyway).

If instead $d < n-1$ we can cover $B_0$ by less than $C M^d$ balls $B$ of radius $C_0$
that are contained in $B(0,M)$, and \eqref{7b17} yields 
$$\frac{\lambda}{c} M^{n-1}\lesssim \, \omega^Y(E \cap B_0)
\leq \sum_B \omega^Y(E \cap B_0) \leq C \frac{\lambda}{c} M^d,$$ 
a contradiction again for $M$ large enough.

This completes our proof of Proposition \ref{t72}. \qed

\begin{remark} \label{r74}
The fact that proving Proposition \ref{t72} was a question of homogeneity
should convince us to try something different
in our Case 1 (elliptic operators): approximate $G^\infty$ with the 
power $D_{\alpha,\mu}(X)^{d+2-n}$, because this is the only one for which \eqref{7b7}
leaves us a fighting chance. (Notice that when $d=n-1$ we keep the power $1$, as in
Theorem \ref{t71}, and this issue does not arise). It turns out such a power is indeed natural in a much more general context, but unfortunately -- perhaps for this reason -- the corresponding results are much more intricate as well. We turn to this question in the next section.
\end{remark}

\section{The condition $\Upsilon_{\rm{flat}}$, higher co-dimension, and more exotic free boundary results. }
\label{S8}

Many known free boundary results rely on some ``zero-level" flatness statement, that is, a degenerate version of the hypothesis entails that the set is a hyperplane. In our context, such a degenerate hypothesis is an equality between the Green function and the distance of our choice. This was not a problem in Section~\ref{S6}: if the Green function is assumed to be the  
distance to a plane, then of course the boundary is a plane simply because the Green function ought to vanish on it. However, when working with $D_{\alpha}$, an analogous statement is not as clear: both the Green function and $D_\alpha$ vanish on the boundary $E$, whatever $E$ is. We were lucky not to have to prove any ``zero-level" flatness in Section~\ref{S7}, by virtue of the results in \cite{HMU}, and, as far as the harmonic measure goes, a straightforward zero-level flatness statement would not even be true -- see Remark~\ref{r71}, but as we pointed out in the same remark, it seems reasonable to expect that such a statement would nonetheless hold for the Green function. Unfortunately, we do not manage to prove it at this point for the desired class of operators at hand, and will simply state it as a hypothesis. 

\begin{definition} \label{d71}
Let the dimension $0 < d < n$, with $d \neq n-2$, $\alpha > 0$,
and a constant coefficient elliptic operator $L_0 = -\dv A_0 \nabla$ be given.  
We say that \ub{the property $\Upsilon_{\rm{flat}}(d, \alpha, L_0)$} is true when 
for each domain $\Omega \subset \R^n$ such that $E = \d \Omega$ is AR 
of dimension $d$ and \eqref{2b4} is satisfied when $d \geq n-1$, the following holds. 
If $\mu$ is an AR measure of dimension $d$ on $E$ and
\begin{equation} \label{7a3}
L_0 D_{\alpha, \mu}^{d+2-n} = 0 \ \text{ on } \Omega,
\end{equation}
then $d$ is an integer, $E$ is a $d$-plane, and $\mu$ is a flat measure.
\end{definition}

\begin{remark} \label{r75}
Our forthcoming results will primarily rely on the case when $L_0=-\Delta$.
In fact, when $L$ is associated to a constant coefficient operator $L_0 \neq - \Delta$,
one could object that our definition of $D_{\alpha,\mu}$ is not as appropriate
as when $L_0 = - \Delta$. That is, maybe we should use different distances like the
$D_\alpha$, but adapted to the operators  $L_0$ that we intend to get at the limit.
This could be more appropriate with respect to the rotation invariance, 
but it looks a little too far-fetched for this paper (and we do not want to have too many definitions). 
\end{remark}

Let us now discuss some particular cases and equivalent reformulations of the condition $\Upsilon_{\rm{flat}}(d, \alpha, L_0)$.
Even when $d=n-1$ and $L_0=-\Delta$, we do not know whether $\Upsilon_{\rm{flat}}(n-1, \alpha, L_0)$ is valid. When $d\in (n-2, n)\setminus \{n-1\}$, the power in \eqref{7a3} is exactly the one that naturally arises in the proof of Proposition~\ref{t72}, see Remark~\ref{r74}. 
We could try to define an analogue of $\Upsilon_{\rm{flat}}(n-2, \alpha, L_0)$ (notice that the case $d=n-2$, as is, does not make much sense), but at this time we will not try to do so.
When 
$d<n-2$, the exponent $d+2-n$ in Definition~\ref{d71} is negative, 
and $D_{\alpha,\mu}^{d+2-n}$ becomes a function with a singularity near $E$.

To show that Definition~\ref{d71} is reasonable we have to verify, at the very least, that \eqref{7a3} is satisfied if $d$ is an integer, $E$ is a $d$-plane, and $\mu$ is a flat measure. And indeed, if $d$ is an integer and $\mu = c \H^d_{\vert P}$ is a flat measure, the scale invariance yields
$R_{\alpha,\mu}(X) = C\dist(X,E)^{-1/\alpha}$, then $D_{\alpha,\mu}(X) = C\dist(X,E)$,
and $D_{\alpha,\mu}^{d+2-n} = C\dist(X,E)^{d+2-n}$. 
When $d=n-1$ and $\Omega$ is a half-space, 
$D_{\alpha, \mu}^{d+2-n}=C t$ ($t$ denoting the coordinate perpendicular to the boundary),
which 
is a solution to any constant coefficient elliptic equation, as desired. 
When $d < n-2$, we denote by $t$
the projection of $X \in \R^n$ on $P^\perp \simeq \R^{n-d}$, and we see that 
$D_{\alpha,\mu}^{d+2-n} = C |t|^{d+2-n}$ is a function of $t$ only.
We can recognize a multiple of the Green function in $\R^{n-d}$. So $D_{\alpha,\mu}^{d+2-n}$
is indeed harmonic in that case. Alternatively, one could make a direct computation to check that $\Delta_{x,t} |t|^{d+2-n} =0$ for $t\in \R^{n-d}\setminus 0$. 
These computations ensure that the condition $\Upsilon_{\rm{flat}}$ is, at least, coherent. Does it have a reasonable chance to be true?

Let us look at an equivalent reformulation of the condition $\Upsilon_{\rm{flat}}$.

\begin{lemma}\label{lU1} Assume that $\Omega\subset \RR^n$ is a domain with an unbounded $d$-Ahlfors regular boundary $E$, $0<d<n$, and $\mu$ is a $d$-AR measure on $E$. Then for any $\alpha>0$
\begin{equation} \label{8a1}
\Delta D_{\alpha, \mu}^{d+2-n} = 0 \, \text{ on $\Omega \,$ if and only if }
L_{_{\alpha,\mu}} D_{\alpha, \mu} = 0 \text{ on $\Omega$,}
\end{equation}
where $L_{_{\alpha,\mu}} $ is given by \eqref{1a2} with $D_\alpha=D_{\alpha, \mu}$.

In particular, for $0 < d < n$, with $d \neq n-2$, and $\alpha > 0$, 
the property  $\Upsilon_{\rm{flat}}(d, \alpha, \Delta)$ is equivalent to the following.
For each domain $\Omega \subset \R^n$ such that $E = \d \Omega$ is AR 
of dimension $d$ and \eqref{2b4} holds if $d \geq n-1$, if $\mu$ is an AR 
measure of dimension $d$ on $E$ and
\begin{equation} \label{7a3-bis}
L_{_{\alpha,\mu}} D_{\alpha, \mu} = 0 \ \text{ on } \Omega,
\end{equation}
then $d$ is an integer, $E$ is a $d$-plane, and $\mu$ is a flat measure.

\end{lemma}

\bp The fact that \eqref{8a1} holds is the result of a direct computation. Notice that all these functions are smooth on $\Omega$, so we can talk about
strong solutions. The verification is easy, because
$$\nabla (D_{\alpha, \mu}^{d+2-n}) = (d+2-n) D_{\alpha, \mu}^{d+1-n} \nabla D_{\alpha, \mu},$$
so $\Delta D_{\alpha, \mu}^{d+2-n} = 0$ if and only if
$\dv D_{\alpha, \mu}^{d+1-n} \nabla D_{\alpha, \mu} = 0$; 
then we just need to compare with \eqref{1a2}. 
Hence, $\Upsilon_{\rm{flat}}(d, \alpha, \Delta)$ says that $D_{\alpha, \mu}$ is only a
solution of $L_{\alpha, \mu}$ in the trivial case of a flat measure, which yields $(1)$ in the statement of the lemma.  
\ep

There is one emblematic case when Lemma \ref{lU1}'s version of 
$\Upsilon_{\rm{flat}}(d, \alpha, \Delta)$ fails miserably;
this is when $\alpha = n-d-2$, the ``magic $\alpha$" case that we discussed in the introduction. In this case, both conditions in \eqref{8a1} hold on any $d$-Ahlfors regular set, which does not even need to have an integer dimension, much less possess any regularity or uniform rectifiability. In the context of this paper this case is certainly not amenable to any free boundary results besides Theorem~\ref{t61}, not only because of the failure of $\Upsilon_{\rm{flat}}$, but because the Green function of $L=-\dv D_{\alpha, \mu}^{-(n-d-1)} \nabla$ is a multiple of $D_{\alpha,  \mu}$ and hence, the harmonic measure is equivalent 
to the Hausdorff measure of the boundary for any $d$-AR set $E$ \cite{DEM}.  It remains to be seen  whether $\Upsilon_{\rm{flat}}(n-2, \alpha, \Delta)$ holds in all other instances.

At this point, let us pass to our main results assuming Condition $\Upsilon_{\rm{flat}}$. As usual, we start with Case 1. 

\begin{theorem} \label{t73} Let $d \in (n-2,n) \sm \{n-1\}$ and $\alpha > 0$ be given, and suppose that the property 
$\Upsilon_{\rm{flat}}(d, \alpha, L_0)$ is true for every constant coefficient elliptic operator 
$L_0 = -\dv A_0 \nabla$. 
Assume that $E$ is Ahlfors regular of dimension $d$, $L$ satisfies \eqref{2b2}, 
\eqref{2b3}, and \eqref{3b3}, and $\Omega$ satisfies \eqref{2b4} if $d > n-1$.
Then the Green function $G^\infty$ for $L$ on $\Omega$ is
not prevalently close to $D_{\alpha,\mu}^{d+2-n}$ for any AR 
measure $\mu$ of dimension $d$ on $E$. 
\end{theorem}

The question is already interesting when $L$ is the Laplacian, and then of course 
Property $\Upsilon_{\rm{flat}}(d, \alpha, \Delta)$ is enough  (see also Remark \ref{r75}).

As we mentioned above,
when $d=n-1$, we do not know whether $\Upsilon_{\rm{flat}}(n-1, \alpha, L_0)$ is true
or not, but we proved Theorem \ref{t71} by different means.

\smallskip
\bp
We shall again prove this by contradiction and compactness. If 
$G^\infty$ is prevalently close to $D_{\alpha,\mu}^{d+2-n}$, this means that 
we can find pairs $(x,r)$ in the analogue $\cG^\ast(\varepsilon, M)$ of 
$\cG_{GD_\alpha}(\varepsilon, M)$ but with the power $d+2-n$, 
with arbitrary values of $\varepsilon$ and $M$. 
We can even choose $(x,r) \in \cG_{cc}(\varepsilon, M)$ too, 
because one of our assumptions, \eqref{3b3}, says that this condition is prevalent. 
We want to show that this is impossible, and in fact we will prove a little more: for each choice of constants in the assumptions of the theorem,
we can find $\varepsilon$ and $M$ such that if $\Omega, E, L$ satisfy the assumptions, 
$\cG^\ast(\varepsilon, M) \cap \cG_{cc}(\varepsilon, M)$ is empty. 

We prove this by contradiction, so we suppose that for each $k \geq 0$, we can find
a counterexample $(\Omega_k, E_k, L_k)$ where this fails for $\varepsilon_k = 2^{-k}$ and 
$M_k = 2^k$, and we want to derive a counterexample. 

We are given a pair
$(x_k, r_k) \in \cG^\ast(\varepsilon_k, M_k) \cap \cG_{cc}(\varepsilon_k, M_k)$,
and by translation and dilation invariance we may assume that $x_k=0$ and $r_k=1$.
Then we proceed as in Theorem \ref{t71}, for instance, and extract a subsequence for which 
$(\Omega_k, E_k, L_k)$ converges to a triple $(\Omega_\infty, E_\infty, L_0)$.
As before, this triple satisfies the assumptions of the theorem (because the $(\Omega_k, E_k, L_k)$ satisfy it uniformly), $L_0$ is a constant coefficient elliptic operator
(because $(0,1) \in \cG_{cc}(\varepsilon_k, M_k)$), and we even know that after 
renormalization the $G_k^\infty$ converge, uniformly on compact subsets of 
$\Omega_\infty$, to a Green function $G$ for $L_0$ on $\Omega_\infty$.
Finally, since $(0,1) \in \cG^\ast(\varepsilon_k, M_k)$, we see that 
$G = D_{\alpha, \mu_\infty}^{d+2-n}$ for some $d$-AR  
measure $\mu_\infty$ that we obtain  as a weak limit of the $\mu_k$ that were used to define 
$\cG^\ast(\varepsilon_k, M_k)$ on $\Omega_k$.

But our assumption $\Upsilon_{\rm{flat}}(d, \alpha, L_0)$ says that this cannot happen
(recall that $d \neq n-1$, so flat measures are not allowed). 
Theorem \ref{t73} follows.
\qed

\ms
Now we pass to Case 2.
As before, we consider the question of whether the Green function for $L_{\alpha,\mu}$
is prevalently close to distance functions, and the most natural ones in this instance 
are clearly the distances $D_{\alpha, \mu}$ associated to the same $\alpha > 0$ and $\mu$ -- cf. \eqref{7a3-bis}.

The natural conjecture is now that $\Upsilon_{\rm{flat}}(d, \alpha, \Delta)$ holds, unless
$d+\alpha = n-2$ or $d=n-2$ (where it is not defined). At this time 
we have some additional information on this question, but no full proof, and for the moment we just state the consequence of a positive result.

\begin{theorem} \label{t81} Let $d \in (0,n) \sm \{n-2\}$ and $\alpha > 0$, $\alpha\neq n-2-d$, be given, and suppose that the property $\Upsilon_{\rm{flat}}(d, \alpha, \Delta)$ holds. 
Let  furthermore $\Omega$ be a domain in $\R^n$, $E = \d \Omega$ be $d$-Ahlfors regular, 
and  assume that \eqref{2b4} holds if $d \geq n-1$.
Let $L = L_{\alpha,\mu} = -\dv A D_{\alpha,\mu}^{d+1-n}\nabla$ be
the degenerate elliptic operator associated to an AR 
measure $\mu$ on $E$ as in \eqref{1a3} and \eqref{1a4}, with the matrix of coefficients $A$ satisfying \eqref{2b2}, \eqref{2b3}, and \eqref{3b3} with $A_0\equiv I$.
If the Green function $G^\infty$ for $L_{\alpha,\mu}$, with pole at $\infty$,
is prevalently close to the distance function $D_{\alpha,\mu}$ 
(with the same $\alpha$ and $\mu$) then
$d$ is an integer and $E$ is uniformly rectifiable.
\end{theorem}

Again $d+\alpha = n-2$ is excluded here, because in this case $D_{\alpha,\mu}$ is the Green function no matter what.

We already know that the result holds 
when $d = n-1$, because in this case $L_{\alpha,\mu} = \Delta$
and we have Theorem \ref{t71}, which fortunately does not need 
$\Upsilon_{\rm{flat}}(d, \alpha, \Delta)$.

\smallskip
\bp
We prove the theorem by compactness as usual, and we start as in Theorem \ref{t71}.
Let $\Omega$, $E$, $L$ satisfy the assumptions. We 
want to prove that $d$ is an integer
and that $E$ is uniformly rectifiable, i.e. that for all $\tau > 0$,
the set $\cG_{ur}(\tau)$ (associated to the integer $d$) is prevalent.
By assumption, $\cG_{GD_{\alpha,\mu}}(\varepsilon, M)$ is prevalent, and we added
$\mu$ in the notation because this time we really need $D_{\alpha,\mu}$ to be defined
in terms of $\mu$. As before, it will be enough to prove that if $\varepsilon$ and $M^{-1}$
are chosen small enough, depending on $n$, $d$, $\tau$, the AR constant for $\mu$, and 
the constant in \eqref{2b4} when $d \geq n-1$, 
\begin{equation} \label{8a3}
\cG_{GD_{\alpha,\mu}}(\varepsilon, M) = \emptyset 
\text{ unless $d$ is an integer}
\end{equation}
and, when $d$ is an integer, 
\begin{equation} \label{8a4}
\cG_{GD_{\alpha,\mu}}(\varepsilon, M) \subset \cG_{ur}(\tau).
\end{equation}
We proceed by contradiction, assume that \eqref{8a3}  or \eqref{8a4} fails for
$\varepsilon_k = 4^{-k}$ and $M_k= 2^{k}$, and pick a counterexample for each $k$.
This means, we find an open set $\Omega_k$ bounded by an AR set $E_k$, 
an AR measure measure $\mu_k$ on $E_k$ (with uniform AR bounds),
and finally a pair $(x_k,r_k) \in E_k \times (0,+\infty)$ such that
the assumptions above hold
(and in particular $(x_k,r_k) \in \cG_{GD_{\alpha,\mu}}(4^{-k}, 2^k)$), 
but not the conclusions.

By translation and dilation invariance, we may assume that $x_k = 0$ and 
$r_k = 1$. We may also extract a subsequence so that $\Omega_k$
converges to a limit $\Omega_\infty$, $E_k$ converges to $E_\infty$, and
$\mu_k$ converges weakly to a limit $\mu_\infty$ on $E_\infty$.
As in the proof of Theorem \ref{l2a4}, $\Omega_\infty$ is bounded by $E_\infty$,
$\Omega_\infty$ satisfies the assumption \eqref{2b4} when needed,
and $\mu_\infty$ is automatically an AR measure on $E_\infty$.

The coefficients of the operator $L_k$  converge in $L^1_{loc} (\Omega_\infty)$
to the coefficients $D_{\alpha,\mu_\infty}^{d+1-n}$ of the operator $L_\infty$
associated to $D_{\alpha,\mu_\infty}$ on $\Omega_\infty$
(see the argument near \eqref{eqConv}: we are now proceeding as in the proof of 
Theorem \ref{t3a1bis}, but with different assumptions). 
In particular, we have \eqref{2a10bis} with $A_\infty=I$, and we are allowed to apply 
Theorem \ref{l2a4}. Hence the Green function $G_k^\infty$ for $L_k$
converges, uniformly on compact subsets of $\Omega_\infty$ to a multiple
of the Green function $G$ for $L_\infty$.

But our assumption that $(0,1) \in \cG_{GD_{\alpha,\mu}}(4^{-k}, 2^k)$ implies that
\begin{equation} \label{8a5}
|D_{\alpha,\mu_k}(X) - c_k G_k^\infty(X)| \leq 2^{-k}
\text{ for } X \in \Omega_k \cap B(0,2^k)
\end{equation}
(compare with \eqref{1a12}), and since $D_{\alpha,\mu_k}$ converges to
$D_{\alpha,\mu_\infty}$, we see that $G = C D_{\alpha,\mu_\infty}$ on $\Omega_\infty$.
In particular $L_\infty D_{\alpha,\mu_\infty} = 0$, and now \eqref{8a1} and our assumption
$\Upsilon_{\rm{flat}}(d, \alpha, \Delta)$ imply that $d$ is an integer, $\mu_\infty$ is a flat measure,
and of course $E_\infty$ is a $d$-plane.

So \eqref{8a3} did not fail, 
and the fact that $E_k$ tends to $E_\infty$
implies that $d_{0,1}(E,E_k)$ tends to $0$, so \eqref{8a4} is also true for $k$ large.
This contradiction completes our proof of \eqref{8a3} and \eqref{8a4},
and Theorem \ref{t81} follows.
\qed

\end{document}